\documentclass[amscd, leqno, amstex]{amsart}
\usepackage{amscd}
\usepackage{amsmath}
\usepackage{amssymb}
\theoremstyle{plain}
\newtheorem{theorem}{Theorem}[section]
\newtheorem{proposition}[theorem]{Proposition}
\newtheorem{lemma}[theorem]{Lemma}

\newtheorem{corollary}[theorem]{Corollary}
\theoremstyle{definition}
\newtheorem{definition}[theorem]{Definition}

\theoremstyle{remark}

\newtheorem{remark}[theorem]{Remark}

\newenvironment{pf}{\begin{proof}}{\end{proof}}
\makeatletter

\@addtoreset{equation}{section}
\makeatother

\begin{document}

\title
{Degenerations of Calabi-Yau threefolds and BCOV invariants}
\author{Ken-Ichi Yoshikawa}
\address{Department of Mathematics, 
Faculty of Science,
Kyoto University,
Kyoto 606-8502, JAPAN}
\email{yosikawa@math.kyoto-u.ac.jp}

\begin{abstract}
In \cite{BCOV93}, \cite{BCOV94}, by expressing the physical quantity $F_{1}$ in two distinct ways,
Bershadsky-Cecotti-Ooguri-Vafa discovered a remarkable equivalence between Ray-Singer analytic torsion and 
elliptic instanton numbers for Calabi-Yau threefolds. 
After their discovery, in \cite{FLY08}, a holomorphic torsion invariant for Calabi-Yau threefolds corresponding to $F_{1}$, called BCOV invariant, was constructed.
In this article, we study the asymptotic behavior of BCOV invariants for algebraic one-parameter degenerations of Calabi-Yau threefolds.
We prove the rationality of the coefficient of logarithmic divergence and give its geometric expression by using a semi-stable reduction of the given family.
\end{abstract}

\maketitle

\tableofcontents

\section*{Introduction}\label{sect:0}
\par
In \cite{BCOV93}, \cite{BCOV94}, by expressing the physical quantity $F_{1}$ in two distinct ways,
Bershadsky-Cecotti-Ooguri-Vafa discovered a remarkable equivalence between Ray-Singer analytic torsion and 
elliptic instanton numbers for Calabi-Yau threefolds. 
After their discovery, in \cite{FLY08}, a holomorphic torsion invariant for Calabi-Yau threefolds corresponding to $F_{1}$, called {\em BCOV invariant}, 
was constructed.
Because of its invariance property, BCOV invariant gives rise to a function $\tau_{\rm BCOV}$ on the moduli space of Calabi-Yau threefolds.
In physics literatures, $-\log\tau_{\rm BCOV}$ is denoted by $F_{1}$.
The prediction of Bershadsky-Cecotti-Ooguri-Vafa concerning the equivalence of holomorphic torsion and elliptic instanton numbers 
for Calabi-Yau threefolds can be stated as follows: 
$\tau_{\rm BCOV}$ admits an explicit infinite product expression of Borcherds type near the large complex structure limit point of the compactified moduli space
of Calabi-Yau threefolds, and the exponents of the infinite product are given by explicit linear combinations of rational and elliptic instanton numbers 
of  the mirror Calabi-Yau threefold corresponding to the large complex structure limit point.
\par
As was done in \cite{BCOV93}, \cite{BCOV94}, a possible first step towards the conjecture of Bershadsky-Cecotti-Ooguri-Vafa
is to determine the section of certain holomorphic line bundle on the moduli space corresponding to $\tau_{\rm BCOV}$.
Thanks to the curvature theorem of Bismut-Gillet-Soul\'e \cite{BGS88}, the complex Hessian $dd^{c}\log\tau_{\rm BCOV}$ is expressed as an explicit linear combination
of the Weil-Petersson form and its Ricci-form on the moduli space \cite{BCOV93}, \cite{BCOV94}, \cite{FangLu03}, \cite{FLY08}.
However, since the moduli space of Calabi-Yau threefolds are non-compact in general, the complex Hessian $dd^{c}\log\tau_{\rm BCOV}$ does not determine
uniquely its potential and hence the corresponding holomorphic section. 
To determine its potential up to a constant, $dd^{c}\log\tau_{\rm BCOV}$ must be determined as a current on some compactified moduli space. 
In this way, we are led to the following two problems: 
one is to understand the behaviors of Weil-Petersson and its Ricci forms as well as their potentials
near the boundary locus of the compactified moduli space;
the other is to understand the behavior of $\tau_{\rm BCOV}$ near the boundary locus of the compactified moduli space. 
We refer to \cite{Lu01}, \cite{LuSun04}, \cite{FLY08}, \cite{LuDouglas13} for the first problem.
In this article, we focus on the second problem.
\par
In this direction,
in \cite{FLY08}, the following results were obtained as an application of the theory of Quillen metrics \cite{BismutLebeau91}, \cite{Bismut97}, \cite{Yoshikawa07}:
$\log\tau_{\rm BCOV}$ always has logarithmic singularity for arbitrary algebraic one-parameter degenerations of Calabi-Yau threefolds and
the logarithmic singularity of $\log\tau_{\rm BCOV}$ is determined for smoothings of Calabi-Yau varieties with at most one ordinary double point
under an additional assumption of the dimension of moduli space.
These results, together with the formula for $dd^{c}\log\tau_{\rm BCOV}$ and the known boundary behaviors
of Weil-Petersson and its Ricci forms, are sufficient to determine $\tau_{\rm BCOV}$ for quintic mirror threefolds \cite{FLY08}. 
However, the results in \cite{FLY08} concerning the singularity of $\tau_{\rm BCOV}$ are not sufficient to determine an explicit formula for $\tau_{\rm BCOV}$ 
for wider classes of Calabi-Yau threefolds, e.g. Calabi-Yau threefolds of Borcea-Voisin. 
For this reason, it is strongly desired to improve the above results in \cite{FLY08}.
The purpose of the present article is to give such improvements. Let us explain our main results.
\par
Let $f\colon{\mathcal X}\to C$ be a surjective morphism from an irreducible projective fourfold ${\mathcal X}$ to a compact Riemann surface $C$.
Assume that there exists a finite subset $\Delta_{f}\subset C$ such that 
$f|_{C\setminus\Delta_{f}}\colon{\mathcal X}|_{C\setminus\Delta_{f}}\to C\setminus\Delta_{f}$ is a smooth morphism 
and such that $X_{t}=f^{-1}(t)$ is a Calabi-Yau threefold for all $t\in C\setminus\Delta_{f}$.

\begin{theorem}
\label{main:theorem:1}
For every $0\in\Delta_{f}$, there exists $\alpha=\alpha_{0}\in{\bf Q}$ such that
$$
\log\tau_{\rm BCOV}(X_{t})=\alpha\,\log|t|^{2}+O\left(\log(-\log|t|)\right)
\qquad
(t\to0),
$$
where $t$ is a local parameter of $C$ centered at $0\in\Delta_{f}$. 
\end{theorem}

We remark that in the corresponding theorem in \cite{FLY08}, the rationality of $\alpha$ was missing.
After Theorem~\ref{main:theorem:1}, a natural question is how the coefficient $\alpha$ is determined by the family $f\colon{\mathcal X}\to C$.
For this, following \cite{FLY08}, we consider its semi-stable reduction \cite{Mumford73}.
Let $g\colon({\mathcal Y},Y_{0})\to(B,0)$ be a semi-stable reduction of $f\colon({\mathcal X},X_{0})\to(C,0)$.
By definition, ${\mathcal Y}$ is a smooth projective fourfold, $(B,0)$ is a pointed compact Riemann surface and there is a surjective morphism of 
pointed compact Riemann surfaces $\phi\colon(B,0)\to(C,0)$ such that $Y_{0}=g^{-1}(0)$ is a {\em reduced} normal crossing divisor of ${\mathcal Y}$ 
and such that ${\mathcal Y}\setminus Y_{0}\cong({\mathcal X}\setminus X_{0})\times_{C\setminus\{0\}}(B\setminus\{0\})$.
\par
By choosing a small neighborhood $V$ of $0$ in $B$, $g^{-1}(V)$ carries a canonical form whose zero divisor is contained in $Y_{0}$.
The zero divisor of any canonical form on $g^{-1}(V)$ with this property is independent of the choice of such canonical form
and is denoted by ${\frak K}_{({\mathcal Y},Y_{0})}$. We call ${\frak K}_{({\mathcal Y},Y_{0})}$ the normalized canonical divisor.
\par
Let $\Omega_{{\mathcal Y}/B}^{1}$ be the sheaf of relative K\"ahler differentials on ${\mathcal Y}$ and 
let $\Omega_{{\mathcal Y}/B}^{1}(\log)$ be its logarithmic version.
Set ${\mathcal Q}=\Omega^{1}_{{\mathcal Y}/B}(\log)/\Omega^{1}_{{\mathcal Y}/B}$.
Every direct image $R^{q}g_{*}{\mathcal Q}|_{V}$ is a finitely generated torsion sheaf on $V$ supported at $0\in B$. We set
$\chi(Rg_{*}{\mathcal Q}|_{V})=\sum_{q\geq0}(-1)^{q}\dim_{\bf C}(R^{q}g_{*}{\mathcal Q})_{0}\in{\bf Z}$.
\par
Let $\varSigma_{g}$ be the critical locus of $g$. 
Let ${\bf P}(T{\mathcal Y})^{\lor}$ be the projective bundle over ${\mathcal Y}$ whose fiber ${\bf P}(T{\mathcal Y})^{\lor}_{y}$ 
is the projective space of hyperplanes of $T_{y}{\mathcal Y}$. 
Then the Gauss map $\mu\colon{\mathcal Y}\setminus\varSigma_{g}\ni y\to[T_{y}Y_{g(y)}]\in{\bf P}(T{\mathcal Y})^{\lor}$
extends to a meromorphic map from ${\mathcal Y}$ to ${\bf P}(T{\mathcal Y})^{\lor}$.
Namely, there exists a blowing-up $\sigma\colon\widetilde{\mathcal Y}\to{\mathcal Y}$ inducing an isomorphism
$\widetilde{\mathcal Y}\setminus\sigma^{-1}(\varSigma_{g})\cong{\mathcal Y}\setminus\varSigma_{g}$ such that
the composite $\widetilde{\mu}=\mu\circ\sigma$ extends to a holomorphic map from $\widetilde{\mathcal Y}$ to ${\bf P}(T{\mathcal Y})^{\lor}$.
Set $\widetilde{g}=g\circ\sigma$. We have a new family of Calabi-Yau threefolds $\widetilde{g}\colon\widetilde{\mathcal Y}\to B$,
whose critical locus $\varSigma_{\widetilde{g}}|_{V}$ defines a divisor of $\widetilde{\mathcal Y}$.
\par
Let $U$ be the universal hyperplane bundle over ${\bf P}(T{\mathcal Y})^{\lor}$ and let $H$ be the universal quotient line bundle over ${\bf P}(T{\mathcal Y})^{\lor}$. 
Following \cite{FLY08}, set
$$
a_{p}(g,\varSigma_{g})
=
\sum_{j=0}^{p}(-1)^{p-j}\,
\int_{{\rm Exc}(\sigma)}
\widetilde{\mu}^{*}
\left\{
{\rm Td}(U)\,
\frac{{\rm Td}(c_{1}(H))-e^{-(p-j)c_{1}(H)}}{c_{1}(H)}
\right\}\,
\sigma^{*}{\rm ch}(\Omega^{j}_{\mathcal Y}),
$$
where ${\rm Exc}(\sigma)$ is the exceptional divisor of $\sigma\colon\widetilde{\mathcal Y}\to{\mathcal Y}$, ${\rm Td}(\cdot)$ is the Todd genus, 
and $\Omega_{\mathcal Y}^{j}$ is the holomorphic vector bundle of holomorphic $j$-forms on ${\mathcal Y}$.
Define 
$$
\rho(g,\varSigma_{g})
=
-3a_{0}(g,\varSigma_{g})+2a_{1}(g,\varSigma_{g})-\chi(Rg_{*}{\mathcal Q}|_{V})
+
\frac{1}{12}\int_{\varSigma_{\widetilde{g}}|_{V}}\widetilde{\mu}^{*}c_{3}(U)
\in
{\bf Q},
$$
$$
\kappa(g,\varSigma_{g},{\frak K}_{({\mathcal Y},Y_{0})})
=
\int_{\sigma^{*}{\frak K}_{({\mathcal Y},Y_{0})}}\widetilde{\mu}^{*}c_{3}(U)
\in
{\bf Z}.
$$

\begin{theorem}
\label{main:theorem:2}
The rational number $\alpha$ in Theorem~\ref{main:theorem:1} is given by 
$$
\alpha
=
\frac{1}{\deg\{\phi\colon(B,0)\to(C,0)\}}
\left\{
\rho(g,\varSigma_{g})
-
\frac{1}{12}\kappa(g,\varSigma_{g},{\frak K}_{({\mathcal Y},Y_{0})})
\right\}.
$$
\end{theorem}

Since every algebraic one-parameter degeneration of Calabi-Yau threefolds admits a semi-stable reduction \cite{Mumford73},
in principle, one can compute the singularity of $\tau_{\rm BCOV}$ for those degenerations by Theorems~\ref{main:theorem:1} and \ref{main:theorem:2},
once one knows their semi-stable reductions. In this sense, the problem of understanding the singularity of $\tau_{\rm BCOV}$ is reduced to
the algebro-geometric problem of classifying possible semi-stable degenerations of Calabi-Yau threefolds. 
\par
As an application of Theorems~\ref{main:theorem:1} and \ref{main:theorem:2}, we shall prove certain locality of the singularity of $\tau_{\rm BCOV}$.
Namely, under some additional assumptions about the family $f\colon{\mathcal X}\to C$ (cf. Section~\ref{sect:4} for the required conditions),
the coefficient $\alpha$ in Theorem~\ref{main:theorem:1} depends only on the function germ of $f$ around the critical locus $\varSigma_{f}$.
(See Theorem~\ref{thm:locality:BCOV:invariant} for the precise statement.)
In some cases, this locality is quite powerful, because we have only to compute one particular example to determine the singularity of $\tau_{\rm BCOV}$. 
This locality result plays a crucial role to determine the BCOV invariant for Borcea-Voisin threefolds \cite{Yoshikawa14}.
\par
The strategy to the proof of Theorems~\ref{main:theorem:1} and \ref{main:theorem:2} is quite parallel to that of \cite[Th.\,9.1]{FLY08}.
In \cite{FLY08}, it was proved that the $L^{2}$-metric on the line bundle $\det R^{q}g_{*}\Omega_{{\mathcal Y}/B}^{p}(\log)$ has 
at most an algebraic singularity at the discriminant locus when $p+q=3$. In this article, we shall improve this estimate. 
Namely, under the assumption of semi-stability, the $L^{2}$-metric on $\det R^{q}g_{*}\Omega_{{\mathcal Y}/B}^{p}(\log)$ has at most a logarithmic singularity,
which enables us to determine various inexplicit constants in \cite[\S9]{FLY08} and hence $\alpha$ in Theorem~\ref{main:theorem:1}.
\par
This article is organized as follows.
In Section~\ref{sect:1}, we recall the construction of BCOV invariants.
In Section~\ref{sect:2}, we study the asymptotic behavior of the $L^{2}$-metric on $\det R^{q}g_{*}\Omega_{{\mathcal Y}/B}^{p}(\log)$ 
and prove the key fact that it has at most a mild singularity when $p+q=3$.
In Section~\ref{sect:3}, we prove Theorems~\ref{main:theorem:1} and \ref{main:theorem:2}.
In Section~\ref{sect:4}, we prove the locality of the singularity of $\tau_{\rm BCOV}$.
In Section~\ref{sect:5}, we determine the singularity of $\tau_{\rm BCOV}$ for general one-parameter smoothings of Calabi-Yau varieties with at most
ordinary double points.
\par
{\bf Acknowledgements }
The author thanks Professor Y. Namikawa for helpful discussions about Kodaira-Spencer maps.
The author is partially supported by JSPS Grants-in-Aid (B) 23340017, (A) 22244003, (S) 22224001, (S) 25220701.

\section
{BCOV invariants}
\label{sect:1}
\par

\subsection
{Analytic torsion and Quillen metrics}
\label{sect:1.1}
\par
Let $(M,g)$ be a compact K\"ahler manifold of dimension $d$ with K\"ahler form $\omega$.
Let $\square_{p,q}=(\bar{\partial}+\bar{\partial}^{*})^{2}$ be the Hodge-Kodaira Laplacian acting on $C^{\infty}$ $(p,q)$-forms on $M$
or equivalently $(0,q)$-forms on $M$ with values in $\Omega_{M}^{p}$, where $\Omega_{M}^{1}$ is the holomorphic cotangent bundle of $M$ and
$\Omega_{M}^{p}:=\Lambda^{p}\Omega_{M}^{1}$.
Let $\sigma(\square_{p,q})\subset{\bf R}_{\geq0}$ be the set of eigenvalues of $\square_{p,q}$. The spectral zeta function of $\square_{p,q}$
is defined as 
$$
\zeta_{p,q}(s):=\sum_{\lambda\in\sigma(\square_{p,q})\setminus\{0\}}\lambda^{-s}\dim E(\lambda,\square_{p,q}),
$$
where $E(\lambda,\square_{p,q})$ is the eigenspace of $\square_{p,q}$ corresponding to the eigenvalue $\lambda$.
Then $\zeta_{p,q}(s)$ converges on the half-plane $\{s\in{\bf C};\,\Re s>\dim M\}$, extends to a meromorphic function on ${\bf C}$,
and is holomorphic at $s=0$. By Ray-Singer \cite{RaySinger73}, the {\em analytic torsion} of $(M,\Omega^{p}_{M})$ is the real number defined as
$$
\tau(M,\Omega_{M}^{p})
:=
\exp\{-\sum_{q\geq0}(-1)^{q}q\,\zeta'_{p,q}(0)\}.
$$
Obviously, $\tau(M,\Omega_{M}^{p})$ depends not only on the complex structure of $M$ but also on the metric $g$.
When we emphasis the dependence of analytic torsion on the metric, we write $\tau(M,\Omega_{M}^{p},g)$.
\par
In \cite{BCOV94}, Bershadsky-Cecotti-Ooguri-Vafa introduced the following combination of analytic torsions.

\begin{definition}
\label{def:BCOV:torsion}
The {\em BCOV torsion} of $(M,g)$ is the real number defined as
$$
T_{\rm BCOV}(M,g)
:=
\prod_{q\geq0}\tau(M,\Omega_{M}^{p})^{(-1)^{p}}
=
\exp\{-\sum_{p,q\geq0}(-1)^{p+q}pq\,\zeta'_{p,q}(0)\}.
$$
\end{definition}

If $\gamma$ is the K\"ahler form of $g$, then we often write $T_{\rm BCOV}(M,\gamma)$ for $T_{\rm BCOV}(M,g)$.
In general, $T_{\rm BCOV}(M,g)$ does depend on the choice of K\"ahler metric $g$ and hence is not a holomorphic invariant of $M$.
When $M$ is a Calabi-Yau threefold, it is possible to construct a holomorphic invariant of $M$ from $T_{\rm BCOV}(M,g)$ by multiplying a correction factor.
Following \cite{FLY08}, let us recall the construction of this invariant.

\subsection
{Calabi-Yau threefolds and BCOV invariants}
\label{sect:1.2}
\par
A compact connected K\"ahler manifold $X$ is {\em Calabi-Yau} if $h^{0,q}(X)=0$ for $0<q<\dim X$ and $K_{X}\cong{\mathcal O}_{X}$,
where $K_{X}$ is the canonical line bundle of $X$.
Our particular interest is the case where $X$ is a threefold. 
Let $X$ be a Calabi-Yau threefold. Let $g=\sum_{i,j}g_{i\bar{j}}\,dz_{i}\otimes d\bar{z}_{j}$ be a K\"ahler metric on $X$ and 
let $\gamma=\gamma_{g}:=\sqrt{-1}\sum_{i,j}g_{i\bar{j}}dz_{i}\wedge d\bar{z}_{j}$ be the corresponding K\"ahler form.
Following the convention in Arakelov geometry, we define
$$
{\rm Vol}(X,\gamma):=\frac{1}{(2\pi)^{3}}\int_{X}\frac{\gamma^{3}}{3!}.
$$
The covolume of $H^{2}(X,{\bf Z})_{\rm free}:=H^{2}(X,{\bf Z})/{\rm Torsion}$ with respect to $[\gamma]$ is defined as
$$
{\rm Vol}_{L^{2}}(H^{2}(X,{\bf Z}),[\gamma])
:=
\det(\langle{\bf e}_{i},{\bf e}_{j}\rangle_{L^{2},[\gamma]})_{1\leq i,j\leq b_{2}(X)}.
$$
Here $\{{\bf e}_{1},\ldots,{\bf e}_{b_{2}(X)}\}$ is a basis of $H^{2}(X,{\bf Z})_{\rm free}={\rm Im}\{H^{2}(X,{\bf Z})\to H^{2}(X,{\bf R})\}$ over ${\bf Z}$ 
and $\langle\cdot,\cdot\rangle_{L^{2},[\gamma]}$ is the inner product on $H^{2}(X,{\bf R})$ induced by integration of harmonic forms.
Namely, if ${\mathcal H}{\bf e}_{i}$ denotes the harmonic representative of ${\bf e}_{i}\in H^{2}(X,{\bf R})$ with respect to $\gamma$
and if $*$ denotes the Hodge star operator, then
$$
\langle{\bf e}_{i},{\bf e}_{j}\rangle_{L^{2},[\gamma]}
:=
\frac{1}{(2\pi)^{3}}\int_{X}{\mathcal H}{\bf e}_{i}\wedge*({\mathcal H}{\bf e}_{j}).
$$
The covolume ${\rm Vol}_{L^{2}}(H^{2}(X,{\bf Z}),[\gamma])$ is the volume of real torus $H^{2}(X,{\bf R})/H^{2}(X,{\bf Z})_{\rm free}$
with respect to the $L^{2}$-metric $\langle\cdot,\cdot\rangle_{L^{2},[\gamma]}$ on $H^{2}(X,{\bf R})$.
\par
As the correction term to the BCOV torsion $T_{\rm BCOV}(X,\gamma)$, we introduce a Bott-Chern term.

\begin{definition}
\label{def:Bott:Chern:correction}
For a Calabi-Yau threefold $X$ equipped with a K\"ahler form, define
$$
A(X,\gamma)
:=
\exp\left[-\frac{1}{12}\int_{X}
\log\left(
\sqrt{-1}\frac{\eta\wedge\bar{\eta}}{\gamma^{3}/3!}
\frac{{\rm Vol}(X,\gamma)}{\|\eta\|_{L^{2}}^{2}}
\right)\,c_{3}(X,\gamma)\right],
$$
where $c_{3}(X,\gamma)$ is the top Chern form of $(X,\gamma)$, $\eta\in H^{0}(X,K_{X})\setminus\{0\}$ is a nowhere vanishing canonical form
on $X$ and $\|\eta\|_{L^{2}}$ is its $L^{2}$-norm, i.e.,
$$
\|\eta\|_{L^{2}}
:=
\frac{1}{(2\pi)^{3}}\int_{X}\sqrt{-1}\eta\wedge\overline{\eta}.
$$
\end{definition}

Obviously, $A(X,\gamma)$ is independent of the choice of $\eta\in H^{0}(X,K_{X})\setminus\{0\}$.
We remark that our definition of $A(X,\gamma)$ differs from the one in \cite[Def.\,4.1]{FLY08} by the factor ${\rm Vol}(X,\gamma)^{\chi(X)/12}$,
where $\chi(X)$ denotes the topological Euler number of $X$.
Notice that $A(X,\gamma)=1$ if $\gamma$ is Ricci-flat.

\begin{definition}
\label{def:BCOV:invariant}
The {\em BCOV invariant} of $X$ is the real number defined as
$$
\tau_{\rm BCOV}(X)
:=
{\rm Vol}(X,\gamma)^{-3+\frac{\chi(X)}{12}}{\rm Vol}_{L^{2}}(H^{2}(X,{\bf Z}),[\gamma])^{-1}T_{\rm BCOV}(X,\gamma)\,A(X,\gamma).
$$
\end{definition}

As an application of the curvature formula for Quillen metrics \cite[Th.\,0.1]{BGS88}, we get the invariance property of $\tau_{\rm BCOV}(X)$ in \cite{FLY08}.

\begin{theorem}
\label{thm:FLY}
For a Calabi-Yau threefold $X$, $\tau_{\rm BCOV}(X)$ is independent of the choice of a K\"ahler form on $X$.
\end{theorem}

\begin{pf}
See \cite[Th.\,4.16]{FLY08}.
\end{pf}

After Theorem~\ref{thm:FLY}, we regard $\tau_{\rm BCOV}$ as a function on the moduli space of Calabi-Yau threefolds. 
In this article, we study the behavior of $\tau_{\rm BCOV}$ for algebraic one-parameter families of Calabi-Yau threefolds and 
improve some results in \cite[\S9]{FLY08}.

\section
{Asymptotic behavior of the $L^{2}$-metric on Hodge bundle}
\label{sect:2}
\par
Let $\varDelta:=\{z\in{\bf C};\,|z|<1\}$ be the unit disc and let $\varDelta^{*}:=\varDelta\setminus\{0\}$ be the unit punctured disc. 
Let ${\frak H}:=\{z\in{\bf C};\,\Im z>0\}$ be the complex upper half-plane. We regard ${\frak H}$ as the universal covering of $\varDelta^{*}$
by the map $\varpi\colon{\frak H}\ni z\to\exp(2\pi iz)\in\varDelta^{*}$.
\par
Let $f\colon{\mathcal Z}\to\varDelta$ be a proper surjective holomorphic map from a smooth complex manifold of dimension $n+1$.
We set $Z_{t}:=f^{-1}(t)$ for $t\in\varDelta$.
If $Z_{t}$ is smooth for all $t\in\varDelta^{*}$ and if $Z_{0}$ is a reduced normal crossing divisor of ${\mathcal Z}$,
then the family $f\colon{\mathcal Z}\to\varDelta$ is called a {\em semi-stable degeneration} of relative dimension $n$.
\par
Let $f\colon{\mathcal Z}\to\varDelta$ be a semi-stable degeneration of relative dimension $n$. 
Set $f^{o}:=f|_{\varDelta^{*}}$ and ${\mathcal Z}^{o}:={\mathcal Z}\setminus Z_{0}$.
Let ${\mathcal L}$ be an ample line bundle on ${\mathcal Z}$.
We consider the cohomology of middle degree and set $\ell:=\dim H^{n}(Z_{t},{\bf C})$ for $t\not=0$. 
Assume that 
\begin{center}
{\em $H^{n}(Z_{t},{\bf C})$ consists of primitive cohomology classes with respect to $c_{1}({\mathcal L})|_{Z_{t}}$.}
\end{center}
By the primitivity, each component $H^{p,q}(Z_{t})$, $p+q=n$, carries the $L^{2}$-inner product
\begin{equation}
\label{eqn:L2:metric:middle:degree:primitive}
(u,v)_{L^{2},t}:=(\sqrt{-1})^{p-q}(-1)^{\frac{n(n-1)}{2}}\int_{Z_{t}}u\wedge\overline{v}.
\end{equation}
In particular, $H^{n}(Z_{t},{\bf C})$ is endowed with the $L^{2}$-Hermitian structure, which is independent of the choice of polarization.

\subsection
{$L^{2}$-length of flat section}
\label{sect:2.1}
\par
The holomorphic vector bundle $R^{n}f^{o}_{*}{\bf C}\otimes{\mathcal O}_{\varDelta^{*}}$ is endowed with the Gauss-Manin connection.
Fix a reference point $t_{0}\in\varDelta$ and set $V:=H^{n}(Z_{t_{0}},{\bf C})$.
Fix a basis $\{v_{1},\ldots,v_{\ell}\}$ of $V$, which is unitary with respect to the $L^{2}$-inner product at $t=t_{0}$. 
Since ${\frak H}$ is simply connected, $\varpi^{*}(R^{n}f^{o}_{*}{\bf C})$ is a trivial local system of rank $\ell$ over ${\frak H}$.
By fixing a point $z_{0}\in{\frak H}$ with $t_{0}=\varpi(z_{0})$, each $v_{i}\in\varpi^{*}(R^{n}f^{o}_{*}{\bf C})|_{z_{0}}$ extends uniquely to 
a flat section ${\bf v}_{i}$ of $\varpi^{*}(R^{n}f^{o}_{*}{\bf C})=V\times{\frak H}$ with respect to the Gauss-Manin connection
such that ${\bf v}_{i}(z_{0})=v_{i}$. We regard ${\bf v}_{i}(z)$ as a $V$-valued holomorphic function on ${\frak H}$.
Namely, ${\bf v}_{i}(z)\in{\mathcal O}({\frak H})\otimes_{\bf C}V$.
\par
On $\varpi^{*}(R^{n}f^{o}_{*}{\bf C})=V\times{\frak H}$, the generator of $\pi_{1}(\varDelta^{*})={\bf Z}$ acts as the Picard-Lefschetz transformation:
There exists $T\in{\rm Aut}(V)$ such that for all $z\in{\frak H}$,
$$
{\bf v}_{i}(z+1)=T{\bf v}_{i}(z)
\qquad
(i=1,\ldots,\ell).
$$
Since ${\bf v}_{1}(z),\ldots,{\bf v}_{\ell}(z)$ are flat with respect to the Gauss-Mannin connection, there exist constants $t_{ij}$, $1\leq i,j\leq\ell$,
such that $T{\bf v}_{i}(z)=\sum_{j}t_{ij}{\bf v}_{j}(z)$. Since $f\colon{\mathcal Z}\to\varDelta$ is a semi-stable degeneration, $T=(t_{ij})$ is unipotent.
\par
Let $(\cdot,\cdot)_{L^{2},z}$ be the $L^{2}$-inner product on $H^{n}(Z_{\varpi(z)},{\bf C})$.
Under the identification of $H^{n}(Z_{\varpi(z)},{\bf C})$ with $V=H^{n}(Z_{\varpi(z_{0})},{\bf C})$ via the Gauss-Manin connection, 
$(\cdot,\cdot)_{L^{2},z}$ is regarded as a family of Hermitian structures on $V$ such that 
\begin{equation}
\label{eqn:transformation:rule:L2:norm}
({\bf v}(z+1),{\bf v}'(z+1))_{L^{2},z+1}=(T{\bf v}(z),T{\bf v}'(z))_{L^{2},z}
\end{equation}
for any flat sections ${\bf v}(z)$, ${\bf v}'(z)\in{\mathcal O}({\frak H})\otimes_{\bf C}V$ with respect to the Gauss-Manin connection.
\par
Let $a<b$ be real numbers with $0<b-a<1$.
By \cite[p.46 Prop.\,25]{Griffiths84}, there exist positive constants $C_{1}=C_{1}(a,b)$, $C_{2}=C_{2}(a,b)$, $k=k(a,b)>0$
such that for all $z\in{\frak H}$ with $a<\Re z< b$, $\Im z\gg0$ and for all ${\bf c}=(c_{1},\ldots,c_{\ell})\in{\bf C}^{\ell}$, one has 
\begin{equation}
\label{eqn:estimate:L2:norm:Griffiths}
C_{1}\,(\Im z)^{-k}\|{\bf c}\|^{2}
\leq
\|\sum_{i=1}^{\ell}c_{i}{\bf v}_{i}(z)\|_{L^{2},z}^{2}
\leq 
C_{2}\,(\Im z)^{k}\|{\bf c}\|^{2},
\end{equation}
where $\|{\bf c}\|^{2}=\sum_{i=1}^{\ell}|c_{i}|^{2}=\|\sum_{i=1}^{\ell}c_{i}{\bf v}_{i}(z)\|_{L^{2},z_{0}}^{2}$.

\begin{remark}
\label{remark:SL2:orbit:theorem}
By the ${\rm SL}_{2}$-orbit theorem of Schmid \cite[Th.\,6.6 and its proof]{Schmid73} (in particular, the last equality of \cite[p.253]{Schmid73}),
there is an integer $\nu({\bf c})\in{\bf Z}_{\geq0}$ depending on ${\bf c}$ such that the ratio 
$\|\sum_{i=1}^{\ell}c_{i}{\bf v}_{i}(z)\|_{L^{2},z}^{2}/(\Im z)^{\nu({\bf c})}$ is bounded from below and above by positive constants
when $a<\Re z< b$, $\Im z\gg0$. Hence we indeed have the following better estimate
\begin{equation}
\label{eqn:estimate:L2:norm:Schmid}
C_{1}\,\|{\bf c}\|^{2}
\leq
\|\sum_{i=1}^{\ell}c_{i}{\bf v}_{i}(z)\|_{L^{2},z}^{2}
\leq 
C_{2}\,(\Im z)^{k}\|{\bf c}\|^{2},
\end{equation}
where $k\in{\bf Z}_{\geq0}$. For our later purpose, the weaker estimate \eqref{eqn:estimate:L2:norm:Griffiths}, whose proof is much easier than
that of the ${\rm SL}_{2}$-orbit theorem, is sufficient.
\end{remark}

\par
Let $H(z)$ be the positive-definite $\ell\times\ell$-Hermitian matrix defined as
$$
H(z)
:=
\begin{pmatrix}
({\bf v}_{1}(z),{\bf v}_{1}(z))_{L^{2},z}&\cdots&({\bf v}_{1}(z),{\bf v}_{\ell}(z))_{L^{2},z}
\\
\vdots&\ddots&\vdots
\\
({\bf v}_{\ell}(z),{\bf v}_{1}(z))_{L^{2},z}&\cdots&({\bf v}_{\ell}(z),{\bf v}_{\ell}(z))_{L^{2},z}
\end{pmatrix}.
$$
Let $\lambda_{\min}(z)$ (resp. $\lambda_{\max}(z)$) be the smallest (resp. largest) eigenvalue of the positive-definite Hermitian matrix $H(z)$.
By \eqref{eqn:estimate:L2:norm:Griffiths}, we get
\begin{equation}
\label{eqn:estimate:eigenvalue:1}
C_{1}\,(\Im z)^{-k}\leq\lambda_{\min}(z)\leq\lambda_{\max}(z)\leq C_{2}\,(\Im z)^{k}.
\end{equation}
For a multi-index $I=\{i_{1}<i_{2}<\cdots<i_{r}\}$ with $|I|:=r\leq\ell$, we define
$$
{\bf v}_{I}(z):={\bf v}_{i_{1}}(z)\wedge{\bf v}_{i_{2}}(z)\wedge\cdots\wedge{\bf v}_{i_{r}}(z)\in\Lambda^{r}V.
$$
For $1\leq r\leq\ell$, let $\Lambda^{r}H(z)$ be the positive-definite $\binom{\ell}{r}\times\binom{\ell}{r}$ Hermitian matrix defined as
$$
\Lambda^{r}H(z)
:=
\left(
({\bf v}_{I}(z),{\bf v}_{J}(z))_{L^{2},z}
\right)_{|I|=|J|=r}.
$$
Here $\Lambda^{r}V$ is equipped with the Hermitian structure induced by $(\cdot,\cdot)_{L^{2},z}$, which is again denoted by the same symbol.
By \eqref{eqn:estimate:eigenvalue:1}, we have the following inequality of positive-definite Hermitian endomorphisms on $\Lambda^{r}V$
\begin{equation}
\label{eqn:estimate:eigenvalue:2}
C_{1}^{r}\,(\Im z)^{-kr}I_{\Lambda^{r}V}
\leq
\lambda_{\min}(z)^{r}I_{\Lambda^{r}V}
\leq
\Lambda^{r}H(z)
\leq
\lambda_{\max}(z)^{r}I_{\Lambda^{r}V}
\leq
C_{2}^{r}\,(\Im z)^{kr}I_{\Lambda^{r}V}.
\end{equation}
By \eqref{eqn:estimate:eigenvalue:2}, for all $z\in{\frak H}$ with $a<\Re z< b$ and $\xi=(\xi_{I})_{|I|=r}\in{\bf C}^{\binom{\ell}{r}}$, we get
\begin{equation}
\label{eqn:estimate:L2:norm:Grassmann}
C_{1}^{r}(\Im z)^{-kr}\|\xi\|^{2}
\leq
\|\sum_{|I|=r}\xi_{I}{\bf v}_{I}(z)\|_{L^{2},z}^{2}
\leq
C_{2}^{r}(\Im z)^{kr}\|\xi\|^{2}.
\end{equation}

\subsection
{$L^{2}$-length of the canonical section associated to flat section}
\label{sect:2.2}
\par
From the flat sections ${\bf v}_{1}(z),\ldots,{\bf v}_{\ell}(z)$, one can construct nowhere vanishing $\pi_{1}(\varDelta^{*})$-invariant holomorphic sections 
${\bf s}_{1},\ldots,{\bf s}_{\ell}$ as follows.
Let $N\in{\rm End}(V)$ be the logarithm of the Picard-Lefschetz transformation $T$. 
Since $T$ is unipotent by our assumption, we have $N=\sum_{l\geq1}(-1)^{l+1}(T-1_{V})^{l}/l$. 
Since the entries of $T\in{\rm Aut}(V)$ with respect to the basis $\{{\bf v}_{1}(z),\ldots,{\bf v}_{\ell}(z)\}$ are constant, so are the entries of $N\in{\rm End}(V)$.
We define 
\begin{equation}
\label{eqn:canonical:extension:Hodge:bundle}
{\bf s}_{i}\left(\exp(2\pi\sqrt{-1}z)\right)
:=
e^{-zN}{\bf v}_{i}(z)
=
\sum_{k\geq0}\frac{(-1)^{k}}{k!}z^{k}N^{k}\,{\bf v}_{i}(z)
\in
{\mathcal O}({\frak H})\otimes_{\bf C}V.
\end{equation}
Since ${\bf v}_{i}(z+1)=e^{N}{\bf v}_{i}(z)$, ${\bf s}_{i}$ is $\pi_{1}(\varDelta^{*})$-invariant and descends to a nowhere vanishing holomorphic section 
of $R^{n}f^{o}_{*}{\bf C}\otimes{\mathcal O}_{\varDelta^{*}}$.  
Since the inner product $({\bf s}_{i},{\bf s}_{j})_{L^{2},z}$ is $\pi_{1}(\varDelta^{*})$-invariant by \eqref{eqn:transformation:rule:L2:norm},
it is denoted by $({\bf s}_{i},{\bf s}_{j})_{L^{2},t}$, where $t=\exp(2\pi\sqrt{-1}z)$.
After Schmid \cite[p.235]{Schmid73}, the {\em canonical extension} of $R^{n}f^{o}_{*}{\bf C}\otimes{\mathcal O}_{\varDelta^{*}}$ to $\varDelta$,
denoted by ${\mathcal H}^{n}$, is defined as the holomorphic vector bundle of rank $\ell$ over $\varDelta$
generated by the frame fields $\{{\bf s}_{1},\ldots,{\bf s}_{\ell}\}$
(see \cite{Steenbrink76}, \cite{Zucker84} for algebro-geometric construction):
$$
{\mathcal H}^{n}:={\mathcal O}_{\varDelta}{\bf s}_{1}\oplus\cdots\oplus{\mathcal O}_{\varDelta}{\bf s}_{\ell}.
$$
\par
Since $N$ is nilpotent and has constant entries with respect to $\{{\bf v}_{1}(z),\ldots,{\bf v}_{\ell}(z)\}$, 
there exists by \eqref{eqn:canonical:extension:Hodge:bundle} polynomials $P_{ij}(z)\in{\bf C}[z]$ such that for all  $z\in{\frak H}$
$$
{\bf s}_{i}(e^{2\pi\sqrt{-1}z})=\sum_{j=1}^{\ell}P_{ij}(z)\,{\bf v}_{j}(z)
\quad
(i=1,\ldots,\ell).
$$
Since $\{{\bf s}_{1}(e^{2\pi\sqrt{-1}z}),\ldots,{\bf s}_{\ell}(e^{2\pi\sqrt{-1}z})\}$ is a basis of $H^{n}(Z_{\varpi(\exp(2\pi\sqrt{-1}z))},{\bf C})$, we get
$$
\det
\left( 
P_{ij}(z)
\right)
\not=0
\qquad
(\forall\,z\in{\frak H}).
$$
\par
As before, for a multi-index $I=\{i_{1}<i_{2}<\cdots<i_{r}\}$, we set
$$
{\bf s}_{I}\left(\exp(2\pi\sqrt{-1}z)\right):={\bf s}_{i_{1}}\left(\exp(2\pi\sqrt{-1}z)\right)\wedge\cdots\wedge{\bf s}_{i_{r}}\left(\exp(2\pi\sqrt{-1}z)\right).
$$
Then we have
$$
\sum_{|I|=r}\xi_{I}\,{\bf s}_{I}(e^{2\pi\sqrt{-1}z})
=
\sum_{|J|=r}
\left(
\sum_{|I|=r}
\xi_{I}
\left|
\begin{matrix}
P_{i_{1}j_{1}}(z)&\cdots&P_{i_{1}j_{r}}(z)
\\
\vdots&\ddots&\vdots
\\
P_{i_{r}j_{1}}(z)&\cdots&P_{i_{r}j_{r}}(z)
\end{matrix}
\right|
\right)
{\bf v}_{J}(z),
$$
where $I=\{i_{1}<\ldots<i_{r}\}$, $J=\{j_{1}<\ldots<j_{r}\}$.
We deduce from \eqref{eqn:estimate:L2:norm:Grassmann} that
for all $z\in{\frak H}$ with $a<\Re z< b$ and $\xi=(\xi_{I})_{|I|=r}\in{\bf C}^{\binom{\ell}{r}}$,
\begin{equation}
\label{eqn:estimate:L2:norm:canonical:extension}
C_{1}^{r}(\Im z)^{-kr}
\leq
\frac{\|\sum_{|I|=r}\xi_{I}\,{\bf s}_{I}(e^{2\pi\sqrt{-1}z})\|_{L^{2},z}^{2}}
{\sqrt{\sum_{|J|=r}
\left|
\sum_{|I|=r}
\xi_{I}
\left|
\begin{matrix}
P_{i_{1}j_{1}}(z)&\cdots&P_{i_{1}j_{r}}(z)
\\
\vdots&\ddots&\vdots
\\
P_{i_{r}j_{1}}(z)&\cdots&P_{i_{r}j_{r}}(z)
\end{matrix}
\right|
\right|^{2}}}
\leq
C_{2}^{r}(\Im z)^{kr}.
\end{equation}
\par
We define an invertible $\binom{\ell}{r}\times\binom{\ell}{r}$-matrix $\Lambda^{r}P(z)$ by
$$
\Lambda^{r}P(z)
:=
\left(
\left|
\begin{matrix}
P_{i_{1}j_{1}}(z)&\cdots&P_{i_{1}j_{r}}(z)
\\
\vdots&\ddots&\vdots
\\
P_{i_{r}j_{1}}(z)&\cdots&P_{i_{r}j_{r}}(z)
\end{matrix}
\right|
\right)_{|I|=|J|=r}
\in GL\left({\bf C}^{\binom{\ell}{r}}\right).
$$
Let $\mu_{\min}(z)$ (resp. $\mu_{\max}(z)$) be the smallest (resp. largest) eigenvalue of the positive-definite $\binom{\ell}{r}\times\binom{\ell}{r}$-matrix 
$G(z):={}^{t}(\Lambda^{r}P(z))\overline{\Lambda^{r}P(z)}$. Then we have
\begin{equation}
\label{eqn:estimate:eigenvalue:G(z)}
\mu_{\min}(z)^{-1}\leq{\rm Tr}\{G(z)^{-1}\},
\qquad
\mu_{\max}(z)\leq{\rm Tr}\,G(z).
\end{equation}
For multi-indices $I$, $J$ with $|I|=|J|=r$, let $r_{IJ}(z)\in{\bf C}(z)$ be the $(I,J)$-entry of the $\binom{\ell}{r}\times\binom{\ell}{r}$-matrix $(\Lambda^{r}P(z))^{-1}$.
By the definition of $G(z)$, we get
$$
{\rm Tr}\{G(z)^{-1}\}=\sum_{|I|=|J|=r}|r_{IJ}(z)|^{2},
$$
$$
{\rm Tr}\,G(z)
=
\sum_{|I|=|J|=r}
\left|
\det
\begin{pmatrix}
P_{i_{1}j_{1}}(z)&\cdots&P_{i_{1}j_{r}}(z)
\\
\vdots&\ddots&\vdots
\\
P_{i_{r}j_{1}}(z)&\cdots&P_{i_{r}j_{r}}(z)
\end{pmatrix}
\right|^{2}.
$$
Since $r_{IJ}(z)$ is a rational function in the variable $z$ with $\det(r_{IJ}(z))\not=0$ in ${\bf C}(z)$
and since $P_{ij}(z)$ is a polynomial in the variable $z$ with $\det(P_{IJ}(z))\not=0$ in ${\bf C}[z]$,
there exist $\mu\in{\bf Z}$, $\nu\in{\bf Z}_{\geq0}$ and constants $C_{3},C_{4}>0$ such that
\begin{equation}
\label{eqn:estimate:trace:G(z)}
{\rm Tr}\{G(z)^{-1}\}\leq C_{3}|z|^{2\mu},
\qquad
{\rm Tr}\,G(z)\leq C_{4}\,|z|^{2\nu}
\end{equation}
for all $z\in{\frak H}$ with $|z|\gg1$.
Since
$$
\sum_{|J|=r}
\left|
\sum_{|I|=r}
\xi_{I}
\left|
\begin{matrix}
P_{i_{1}j_{1}}(z)&\cdots&P_{i_{1}j_{r}}(z)
\\
\vdots&\ddots&\vdots
\\
P_{i_{r}j_{1}}(z)&\cdots&P_{i_{r}j_{r}}(z)
\end{matrix}
\right|
\right|^{2}
=
\|\Lambda^{r}P(z)\xi\|^{2}
=
(G(z)\xi,\xi),
$$
we deduce from \eqref{eqn:estimate:eigenvalue:G(z)}, \eqref{eqn:estimate:trace:G(z)} that for $z\in{\frak H}$ with $|z|\gg1$ and 
$\xi=(\xi_{I})_{|I|=r}\in{\bf C}^{\binom{\ell}{r}}$
\begin{equation}
\label{eqn:estimate:denominator}
C_{3}^{-1}\,\|\xi\|^{2}\,|z|^{-2\mu}
\leq
\sum_{|J|=r}
\left|
\sum_{|I|=r}
\xi_{I}
\left|
\begin{matrix}
P_{i_{1}j_{1}}(z)&\cdots&P_{i_{1}j_{r}}(z)
\\
\vdots&\ddots&\vdots
\\
P_{i_{r}j_{1}}(z)&\cdots&P_{i_{r}j_{r}}(z)
\end{matrix}
\right|
\right|^{2}
\leq
C_{4}\,\|\xi\|^{2}\,|z|^{2\nu}.
\end{equation}
By \eqref{eqn:estimate:L2:norm:canonical:extension}, \eqref{eqn:estimate:denominator}, we get
for all $z\in{\frak H}$ with $a<\Re z< b$, $|z|\gg1$ and $\xi\in{\bf C}^{\binom{\ell}{r}}$
\begin{equation}
\label{eqn:estimate:L2:norm:canonical:extension:2}
C_{5}\,\|\xi\|^{2}\,(\Im z)^{-(kr+2\mu)}
\leq
\|\sum_{|I|=r}\xi_{I}\,{\bf s}_{I}(e^{2\pi\sqrt{-1}z})\|_{L^{2},z}^{2}
\leq
C_{6}\,\|\xi\|^{2}\,(\Im z)^{(kr+2\nu)},
\end{equation}
where $C_{5},C_{6}$ are constants. Here we used the fact that $|z|/\Im z$ is bounded from below and above by positive constants 
on the domain $\{z\in{\frak H};\,a<\Re z<b,\,|z|\gg1\}$.
Since $\|{\bf s}_{I}(e^{2\pi\sqrt{-1}z})\|_{L^{2},z}$ is ${\bf Z}$-invariant, it follows from \eqref{eqn:estimate:L2:norm:canonical:extension:2} that
for $t\in\varDelta^{*}$ with $0<|t|\ll1$,
\begin{equation}
\label{eqn:estimate:L2:norm:canonical:extension:3}
C_{5}\,\|\xi\|^{2}\left(-\log|t|\right)^{-(kr+2\mu)}
\leq
\|\sum_{|I|=r}\xi_{I}\,{\bf s}_{I}(t)\|_{L^{2},t}^{2}
\leq
C_{6}\,\|\xi\|^{2}\,(-\log|t|)^{(kr+2\nu)}.
\end{equation}

\subsection
{$L^{2}$-length of a nowhere vanishing section of $\det{\mathcal F}^{p}$}
\label{sect:2.3}
\par
On ${\mathcal H}^{n}$, we have the Hodge filtration
$$
0\subset{\mathcal F}^{n}\subset{\mathcal F}^{n-1}\subset\cdots\subset{\mathcal F}^{1}\subset{\mathcal F}^{0}={\mathcal H}^{n},
$$
where ${\mathcal F}^{p}$ is a holomorphic subbundle of ${\mathcal H}^{n}$ satisfying
${\mathcal F}^{p}_{t}=\bigoplus_{k\geq p}H^{k,n-k}(Z_{t})$ for $t\in\varDelta^{*}$ and
$$
{\mathcal F}^{p}/{\mathcal F}^{p+1}\cong R^{n-p}f_{*}\Omega^{p}_{{\mathcal Z}/\varDelta}(\log Z_{0}).
$$
Here $\Omega^{p}_{{\mathcal Z}/\varDelta}(\log Z_{0}):=\Lambda^{p}\Omega^{1}_{{\mathcal Z}/\varDelta}(\log Z_{0})$ and
$\Omega^{1}_{{\mathcal Z}/\varDelta}(\log Z_{0}):=\Omega^{1}_{\mathcal Z}(\log Z_{0})/{\mathcal O}_{\mathcal Z}f^{*}(dt/t)$.
On the open subset of ${\mathcal Z}$ on which $f(z)=z_{1}\cdots z_{k}$, $\Omega^{1}_{\mathcal Z}(\log Z_{0})$ is given by
$$
\Omega^{1}_{\mathcal Z}(\log Z_{0})
=
{\mathcal O}_{\mathcal Z}(dz_{1}/z_{1})+\cdots+{\mathcal O}_{\mathcal Z}(dz_{k}/z_{k})
+
{\mathcal O}_{\mathcal Z}dz_{k+1}+\cdots+{\mathcal O}_{\mathcal Z}dz_{4}.
$$
See \cite{Steenbrink76}, \cite{Zucker84} for algebro-geometric account of the Hodge filtrations.
In what follows, we often write $R^{q}f_{*}\Omega^{p}_{{\mathcal Z}/\varDelta}(\log)$ for $R^{q}f_{*}\Omega^{p}_{{\mathcal Z}/\varDelta}(\log Z_{0})$.
\par
Set $\ell_{p}:={\rm rk}\,{\mathcal F}^{p}$. There exist holomorphic sections $\varphi_{1},\ldots,\varphi_{\ell}\in\Gamma(\varDelta,{\mathcal H}^{n})$ with
$$
{\mathcal F}^{p}={\mathcal O}_{\varDelta}\varphi_{1}\oplus\cdots\oplus{\mathcal O}_{\varDelta}\varphi_{\ell_{p}}
\qquad
(p=n,n-1,\ldots,1,0).
$$
Since $\{{\bf s}_{1},\ldots,{\bf s}_{\ell}\}$ is a basis of ${\mathcal H}^{n}$ as an ${\mathcal O}_{\varDelta}$-module, 
there exist holomorphic functions $a_{\alpha i}(t)\in{\mathcal O}(\varDelta)$, $1\leq i,\alpha\leq\ell$ such that
$$
\varphi_{\alpha}(t)=\sum_{i=1}^{\ell}a_{\alpha i}(t)\,{\bf s}_{i}(t)
\quad
(\alpha=1,\ldots,\ell),
\qquad
\det\left(a_{\alpha i}(0)\right)_{1\leq\alpha,i\leq\ell}\not=0.
$$

\begin{proposition}
\label{prop:estimate:L2:norm:determinant}
There exist constants $C\in{\bf R}_{\geq0}$ and $C'\in{\bf R}_{>0}$ such that for all $t\in\varDelta^{*}$ with $|t|\ll1$ and $1\leq m\leq\ell$,
$$
\left|
\log
\left\|
\varphi_{1}(t)\wedge\cdots\wedge\varphi_{m}(t)
\right\|_{L^{2},t}
\right|
\leq
C'+C\log\left(-\log|t|\right).
$$
\end{proposition}

\begin{pf}
Since 
$$
\varphi_{1}(t)\wedge\cdots\wedge\varphi_{m}(t)
=
\sum_{|J|=m}
\left|
\begin{matrix}
a_{1j_{1}}(t)&\cdots&a_{1j_{m}}(t)
\\
\vdots&\ddots&\vdots
\\
a_{mj_{1}}(t)&\cdots&a_{mj_{m}}(t)
\end{matrix}
\right|
{\bf s}_{J}(t),
$$
there exist by \eqref{eqn:estimate:L2:norm:canonical:extension:3} constants $k_{1},k_{2}\geq0$ such that
for all $t\in\varDelta^{*}$ with $0<|t|\ll1$
\begin{equation}
\label{eqn:estimate:L2:norm:determinant}
C_{5}\,\left(-\log|t|\right)^{-k_{1}}
\leq
\frac{\|\varphi_{1}(t)\wedge\cdots\wedge\varphi_{m}(t)\|_{L^{2},t}^{2}}
{\sum_{|J|=m}
\left|
\det
\begin{pmatrix}
a_{1j_{1}}(t)&\cdots&a_{1j_{m}}(t)
\\
\vdots&\ddots&\vdots
\\
a_{mj_{1}}(t)&\cdots&a_{mj_{m}}(t)
\end{pmatrix}
\right|^{2}}
\leq
C_{6}\,\left(-\log|t|\right)^{k_{2}}.
\end{equation}
By \eqref{eqn:estimate:L2:norm:determinant}, it suffices to prove that
\begin{equation}
\label{eqn:nonvanishing}
\sum_{|J|=m}
\left|
\det
\begin{pmatrix}
a_{1j_{1}}(0)&\cdots&a_{1j_{m}}(0)
\\
\vdots&\ddots&\vdots
\\
a_{mj_{1}}(0)&\cdots&a_{mj_{m}}(0)
\end{pmatrix}
\right|^{2}
\not=0.
\end{equation}
\par
Set
$$
A
:=
\begin{pmatrix}
a_{11}(0)&\cdots&a_{1\ell}(0)
\\
\vdots&\ddots&\vdots
\\
a_{\ell1}(0)&\cdots&a_{\ell\ell}(0)
\end{pmatrix}
\in
GL({\bf C}^{\ell}).
$$
Then the endomorphism on $\Lambda^{m}{\bf C}^{\ell}$ induced by $A$ is given by the $\binom{\ell}{m}\times\binom{\ell}{m}$-matrix
$$
\Lambda^{m}A
:=
\left(
\left|
\begin{matrix}
a_{i_{1}j_{1}}(0)&\cdots&a_{i_{1}j_{m}}(0)
\\
\vdots&\ddots&\vdots
\\
a_{i_{m}j_{1}}(0)&\cdots&a_{i_{m}j_{m}}(0)
\end{matrix}
\right|
\right)_{|I|=|J|=m}.
$$
Since $\det A\not=0$, $\Lambda^{m}A\in{\rm End}(\Lambda^{m}{\bf C}^{\ell})$ is invertible. In particular, the row vector of $\Lambda^{m}A$ 
corresponding to the multi-index $I=\{1,2,\ldots,m\}$ is non-zero, which implies \eqref{eqn:nonvanishing}. This proves the result.
\end{pf}

Since $R^{n-p}f_{*}\Omega^{p}_{{\mathcal Z}/\varDelta}(\log)|_{\varDelta^{*}}=R^{n-p}f_{*}\Omega^{p}_{{\mathcal Z}/\varDelta}|_{\varDelta^{*}}$
is identified with the holomorphic vector bundle over $\varDelta^{*}$ with fiber $H^{n-p}(Z_{t},\Omega_{Z_{t}}^{p})$ over $t\in\varDelta^{*}$,
$R^{n-p}f_{*}\Omega^{p}_{{\mathcal Z}/\varDelta}(\log)|_{\varDelta^{*}}$ is equipped with the $L^{2}$-Hermitian metric
by identifying $H^{n-p}(Z_{t},\Omega_{Z_{t}}^{p})$ with the corresponding vector space of harmonic forms of bidegree $(p,n-p)$.
\par
Under the canonical identification $R^{n-p}f_{*}\Omega^{p}_{{\mathcal Z}/\varDelta}(\log)={\mathcal F}^{p}/{\mathcal F}^{p+1}$, 
the $L^{2}$-metric on $R^{n-p}f_{*}\Omega^{p}_{{\mathcal Z}/\varDelta}(\log)|_{\varDelta^{*}}$ is identified with
the quotient metric on ${\mathcal F}^{p}/{\mathcal F}^{p+1}$ induced by the $L^{2}$-metric on ${\mathcal H}^{p}|_{\varDelta^{*}}$.
Hence we have an isometry of holomorphic line bundles equipped with singular Hermitian metrics
\begin{equation}
\label{eqn:identification:determinant:bundles}
\left(
\det R^{n-p}f_{*}\Omega^{p}_{{\mathcal Z}/\varDelta}(\log),\|\cdot\|_{L^{2}}
\right)
\cong
(\det{\mathcal F}^{p},\|\cdot\|_{L^{2}})\otimes(\det{\mathcal F}^{p+1},\|\cdot\|_{L^{2}})^{\lor}.
\end{equation}

\begin{corollary}
\label{cor:estimate:L2:norm:det:cohomology}
Let $e_{1}(t),\ldots,e_{n_{p}}(t)\in\Gamma(\varDelta,R^{n-p}f_{*}\Omega^{p}_{{\mathcal Z}/\varDelta}(\log))$ be a basis of 
the free ${\mathcal O}_{\varDelta}$-module $R^{n-p}f_{*}\Omega^{p}_{{\mathcal Z}/\varDelta}(\log)$, where $n_{p}:=h^{p,n-p}(Z_{t})$, $t\not=0$.
Then there exist constants $C\geq0$ and $C'>0$ such that for all $t\in\varDelta^{*}$ with $|t|\ll1$
$$
\left|
\log
\left\|
e_{1}(t)\wedge\cdots\wedge e_{n_{p}}(t)
\right\|_{L^{2},t}
\right|
\leq
C'+C\log\left(-\log|t|\right).
$$
\end{corollary}

\begin{pf}
There exist nowhere vanishing holomorphic sections ${\bf f}_{p}(t)\in\Gamma(\varDelta,\det{\mathcal F}^{p})$ and 
${\bf f}_{p+1}(t)\in\Gamma(\varDelta,\det{\mathcal F}^{p+1})$ such that
$$
e_{1}(t)\wedge\cdots\wedge e_{n_{p}}(t)={\bf f}_{p}(t)\otimes{\bf f}_{p+1}(t)^{-1}
$$
under the identification \eqref{eqn:identification:determinant:bundles}. Since
$$
\|e_{1}(t)\wedge\cdots\wedge e_{n_{p}}(t)\|_{L^{2},t}
=
\|{\bf f}_{p}(t)\|_{L^{2},t}\cdot\|{\bf f}_{p+1}(t)\|_{L^{2},t}^{-1}
$$
the result follows from Proposition~\ref{prop:estimate:L2:norm:determinant} applied to the nowhere vanishing sections ${\bf f}_{p}(t)$ and ${\bf f}_{p+1}(t)$.
\end{pf}

\section
{Singularity of BCOV invariants for semi-stable degenerations} 
\label{sect:3}
\par
{\bf Set up. }
Throughout this section except Section~\ref{sect:3.9}, we assume the following:
\begin{itemize}
\item[(1)]
There exist a smooth projective fourfold ${\mathcal X}$, a compact Riemann surface $C$, an embedding $\varDelta\subset C$,
and a surjective holomorphic map $f\colon{\mathcal X}\to C$ such that $f\colon{\mathcal X}|_{\varDelta}\to\varDelta$ is a semi-stable degeneration.
In particular, $f^{-1}(0)$ is a reduced normal crossing divisor of ${\mathcal X}$.
\item[(2)]
The regular fibers of $f\colon{\mathcal X}\to C$ are Calabi-Yau threefolds.
\end{itemize}
Set $X_{t}:=f^{-1}(t)$ for $t\in C$. Then $X_{t}$ is a Calabi-Yau threefold for all $t\in\varDelta^{*}$.
In this section, we determine the asymptotic behavior of the function on $\varDelta^{*}$
$$
t\mapsto\log\tau_{\rm BCOV}(X_{t})
\qquad
(t\to0).
$$
\par
Let ${\mathcal L}$ be an ample line bundle on ${\mathcal X}$. For all $t\in\varDelta^{*}$, $c_{1}({\mathcal L}_{t})$ is a K\"ahler class on $X_{t}$,
so that every $H^{p,q}(X_{t},{\bf C})$ is endowed with the $L^{2}$-Hermitian structure with respect to the K\"ahler class $c_{1}({\mathcal L}_{t})$.
When $p+q=3$, $H^{p,q}(X_{t},{\bf C})$ consists of primitive cohomology classes and the $L^{2}$-Hermitian structure on $H^{p,q}(X_{t},{\bf C})$ 
with respect to $c_{1}({\mathcal L}_{t})$ is given by \eqref{eqn:L2:metric:middle:degree:primitive}.
\par
Let $\varSigma_{f}:=\{x\in{\mathcal X};\,df_{x}=0\}$ be the critical locus of $f$ and let
$\Delta_{f}:=f(\varSigma_{f})\subset C$ be the discriminant locus of $f\colon{\mathcal X}\to C$. 
Then $\Omega^{p}_{{\mathcal X}/C}$ is a holomorphic vector bundle over ${\mathcal X}\setminus\varSigma_{f}$,
where $\Omega^{p}_{{\mathcal X}/C}:=\Lambda^{p}\Omega^{1}_{{\mathcal X}/C}$ and $\Omega^{1}_{{\mathcal X}/C}:=\Omega^{1}_{\mathcal X}/f^{*}\Omega^{1}_{C}$.
Since possible extension of $\Omega^{p}_{{\mathcal X}/C}$ to a coherent sheaf on ${\mathcal X}$ is not unique in general,
we regard $\Omega^{p}_{{\mathcal X}/C}$ as a locally free sheaf on ${\mathcal X}\setminus\varSigma_{f}$ rather than a coherent sheaf on ${\mathcal X}$.
\par
Recall that the determinant of the cohomologies of $\Omega^{p}_{{\mathcal X}/C}$ is the holomorphic line bundle on $C\setminus\Delta_{f}$ defined as
$$
\lambda(\Omega^{p}_{{\mathcal X}/C})
:=
\bigotimes_{q\geq0}(\det R^{q}f_{*}\Omega^{p}_{{\mathcal X}/C})^{(-1)^{q}}.
$$
Since $X_{t}$, $t\in C\setminus\Delta_{f}$, is equipped with the K\"ahler class $c_{1}({\mathcal L}_{t})$, $H^{q}(X_{t},\Omega^{p}_{X_{t}})$ is equipped with 
the $L^{2}$ Hermitian metric by identifying it with the corresponding vector space of harmonic forms. 
In this way, for all $p,q\geq0$, $R^{q}f_{*}\Omega^{p}_{{\mathcal X}/C}$ is a holomorphic vector bundle over $C\setminus\Delta_{f}$ equipped with
the $L^{2}$ Hermitian metric. (This Hermitian metric coincides with the one considered in Section~\ref{sect:2.3} when $p+q=3$.)
Hence $\lambda(\Omega^{p}_{{\mathcal X}/C})$ is a holomorphic Hermitian line bundle on $C\setminus\Delta_{f}$ equipped with the $L^{2}$ metric 
$\|\cdot\|_{L^{2}}$. 
\par
We fix a K\"ahler metric $g^{\mathcal X}$ on ${\mathcal X}$ with K\"ahler class $c_{1}({\mathcal L})$.
For $t\in C\setminus\Delta_{f}$, we set $g_{t}:=g^{\mathcal X}|_{X_{t}}$. Then $\{(X_{t},g_{t})\}_{t\in C\setminus\Delta_{f}}$ is a family of compact K\"ahler manifolds
with constant K\"ahler class $c_{1}({\mathcal L}_{t})=c_{1}({\mathcal L})|_{X_{t}}$. 
As in Section~\ref{sect:1.1}, we have analytic torsion $\tau(X_{t},\Omega_{X_{t}}^{p})$ for all $t\in C\setminus\Delta_{f}$.
By Quillen \cite{Quillen85} and Bismut-Gillet-Soul\'e \cite{BGS88}, the {\em Quillen metric} on $\lambda(\Omega^{p}_{{\mathcal X}/C})$ is defined as
$$
\|\cdot\|_{Q}^{2}(t):=\tau(X_{t},\Omega_{X_{t}}^{p})\cdot\|\cdot\|_{L^{2},t}^{2},
\qquad
t\in C\setminus\Delta_{f}.
$$
The curvature and anomaly formulae for Quillen metrics were obtained by Bismut-Gillet-Soul\'e \cite{BGS88}.
As an application of the Bismut-Lebeau embedding formula \cite{BismutLebeau91},
the singularity of Quillen metric as $t\to0\in\Delta_{f}$ was determined in \cite{FLY08}, which we shall recall in Section~\ref{sect:3.4}.

\subsection
{Asymptotic behavior of the $L^{2}$-metric on $\lambda({\mathcal O}_{\mathcal X})|_{\varDelta}$}
\label{sect:3.1}
\par

\begin{proposition}
\label{prop:L2:metric:det:cohomology:q=0}
Let $\varsigma_{0}\in\Gamma(\varDelta,\lambda({\mathcal O}_{\mathcal X})|_{\varDelta})$ be a nowhere vanishing holomorphic section.
Then the following holds as $t\to0$:
$$
\log\|\varsigma_{0}(t)\|_{L^{2}}^{2}=O\left(\log(-\log|t|)\right).
$$
\end{proposition}

\begin{pf}
Since $X_{t}$ is a Calabi-Yau threefold for all $t\in\varDelta^{*}$, we have 
$f_{*}{\mathcal O}_{\mathcal X}={\mathcal O}_{C}\cdot1$ 
and
$R^{3}f_{*}{\mathcal O}_{\mathcal X}|_{\varDelta}={\mathcal O}_{\varDelta}\cdot e$, 
where $e(t)\in\Gamma(\varDelta,R^{3}f_{*}{\mathcal O}_{\mathcal X})$ is a nowhere vanishing section. 
By the definition of the $L^{2}$-metric on $H^{0}(X_{t},{\bf C})$, we see that $\|1\|_{L^{2}}^{2}=\deg(L_{t})/(2\pi)^{3}$ is a constant function on $\varDelta$.
By Corollary~\ref{cor:estimate:L2:norm:det:cohomology} in the case $n=3$, $q=0$, we have
$$
\log\|e(t)\|_{L^{2}}^{2}=O\left(\log(-\log|t|)\right)
\qquad
(t\to0).
$$
Since $\sigma_{0}=1\otimes e^{\lor}$, we get the result.
\end{pf}

\subsection
{Asymptotic behavior of the $L^{2}$-metric on $\lambda(\Omega^{1}_{{\mathcal X}/C})|_{\varDelta}$}
\label{sect:3.2}
\par
We define 
$$
\Omega^{1}_{{\mathcal X}/C}:=\Omega^{1}_{\mathcal X}/f^{*}\Omega^{1}_{C}.
$$
Since $\Omega^{1}_{\mathcal X}\subset\Omega^{1}_{\mathcal X}(\log)$ and $\Omega^{1}_{C}\subset\Omega^{1}_{C}(\log)$,
we have the natural inclusion of sheaves
$$
\Omega^{1}_{{\mathcal X}/C}\subset\Omega^{1}_{{\mathcal X}/C}(\log)
$$
and we set
$$
{\mathcal Q}:=\Omega^{1}_{{\mathcal X}/C}(\log)/\Omega^{1}_{{\mathcal X}/C}.
$$
Then ${\mathcal Q}|_{f^{-1}(\varDelta)}$ is a coherent sheaf on ${\mathcal X}|_{\varDelta}$ supported on ${\rm Sing}\,X_{0}$. 
The short exact sequence of coherent sheaves on ${\mathcal X}$
$$
0\longrightarrow\Omega^{1}_{{\mathcal X}/C}\longrightarrow\Omega^{1}_{{\mathcal X}/C}(\log)\longrightarrow{\mathcal Q}\longrightarrow0
$$
induces the long exact sequence of direct image sheaves on $C$
\begin{equation}
\label{eqn:long:exact:sequence:direct:image}
\to
R^{q-1}f_{*}\Omega^{1}_{{\mathcal X}/C}(\log)
\to 
R^{q-1}f_{*}{\mathcal Q}
\to 
R^{q}f_{*}\Omega^{1}_{{\mathcal X}/C}
\to
R^{q}f_{*}\Omega^{1}_{{\mathcal X}/C}(\log)
\to 
R^{q}f_{*}{\mathcal Q}
\to
\end{equation}
Following \cite[proof of Prop.\,9.5]{FLY08}, set
$$
M_{q}:=(R^{q}f_{*}\Omega^{1}_{{\mathcal X}/C})_{\rm tors}|_{\varDelta},
\qquad
N_{q}:=R^{q}f_{*}\Omega^{1}_{{\mathcal X}/C}(\log)/R^{q}f_{*}\Omega^{1}_{{\mathcal X}/C}|_{\varDelta}.
$$
Since $R^{q}f_{*}\Omega^{1}_{{\mathcal X}/C}(\log)$ is a locally free sheaf on $C$,
we deduce from \eqref{eqn:long:exact:sequence:direct:image} the isomorphism of ${\mathcal O}_{\varDelta}$-modules
\begin{equation}
\label{eqn:direct:image:Q}
R^{q}f_{*}{\mathcal Q}|_{\varDelta}
\cong
N_{q}
\oplus
M_{q+1}.
\end{equation}
Set
$$
\chi(Rf_{*}{\mathcal Q}|_{\varDelta}):=\sum_{q}(-1)^{q}\dim_{\bf C}(R^{q}f_{*}{\mathcal Q})_{0}.
$$
By \eqref{eqn:direct:image:Q}, we get
\begin{equation}
\label{eqn:formula:Euler:characteristic}
\chi(Rf_{*}{\mathcal Q}|_{\varDelta})
=
\sum_{q}(-1)^{q}(\dim_{\bf C}(N_{q})_{0}-\dim_{\bf C}(M_{q})_{0}).
\end{equation}
Since ${\mathcal Q}|_{f^{-1}(\varDelta)}$ depends only on the function germ of $f$ near $\varSigma_{f}$, so is $\chi(Rf_{*}{\mathcal Q}|_{\varDelta})$.

\begin{proposition}
\label{prop:L2:metric:det:cohomology:q=1}
Let $\varsigma_{1}\in\Gamma(\varDelta,\lambda(\Omega^{1}_{{\mathcal X}/C}))$ be a nowhere vanishing holomorphic section.
Then the following holds as $t\to0$:
$$
\log
\left\|
\varsigma_{1}(t)
\right\|_{\lambda(\Omega^{1}_{{\mathcal X}/C}),L^{2}}^{2}
=
\chi(Rf_{*}{\mathcal Q}|_{\varDelta})\,\log|t|^{2}+O\left(\log(-\log|t|)\right).
$$
\end{proposition}

\begin{pf}
Let $e_{1}(t),\ldots, e_{h^{1,q}}(t)\in\Gamma(\varDelta,R^{q}f_{*}\Omega^{1}_{{\mathcal X}/C}(\log))$ be a basis of
$R^{q}f_{*}\Omega^{1}_{{\mathcal X}/C}(\log)|_{\varDelta}$ as a free ${\mathcal O}_{\varDelta}$-module.
By \cite[Prop.\,9.4]{FLY08}, there exists $\delta_{q}\in{\bf R}$ such that 
$$
\log\|e_{1}(t)\wedge\cdots\wedge e_{h^{1,q}}(t)\|_{L^{2}}^{2}=\delta_{q}\,\log|t|^{2}+O\left(\log(-\log|t|)\right)
\qquad
(t\to0).
$$
It follows from \cite[Eq.\,(9.15)]{FLY08} and \eqref{eqn:formula:Euler:characteristic} that as $t\to0$
\begin{equation}
\label{eqn:behavior:L2:metric}
\log
\left\|
\varsigma_{1}(t)
\right\|_{\lambda(\Omega^{1}_{{\mathcal X}/C}),L^{2}}^{2}
=
\{\chi(Rf_{*}{\mathcal Q}|_{\varDelta})+\sum_{q}(-1)^{q}\delta_{q}\}\,\log|t|^{2}+O\left(\log(-\log|t|)\right).
\end{equation}
By \eqref{eqn:behavior:L2:metric}, it suffices to prove $\delta_{q}=0$ for all $q\geq0$.
The vanishing $\delta_{0}=\delta_{1}=\delta_{3}=0$ was already proved in \cite[proof of Prop.\,9.4]{FLY08}.
The vanishing $\delta_{2}=0$ follows from Corollary~\ref{cor:estimate:L2:norm:det:cohomology}.
This completes the proof. 
\end{pf}

\subsection
{The K\"ahler extension of $\lambda(\Omega^{p}_{{\mathcal X}/C})|_{\varDelta^{*}}$}
\label{sect:3.3}
\par
Following \cite{FLY08}, we recall an extension of $\lambda(\Omega^{p}_{{\mathcal X}/C})$ from $\varDelta^{*}$ to $\varDelta$, 
which we call the K\"ahler extension and which is distinct from $\lambda(\Omega^{p}_{{\mathcal X}/C}(\log))$ in general.  
For $p\geq 0$, set
$$
\Omega^{p}_{{\mathcal X}/C}:=\Lambda^{p}\Omega^{1}_{{\mathcal X}/C}.
$$
Then $\Omega^{p}_{{\mathcal X}/C}$ is a coherent sheaf on ${\mathcal X}$, which is locally free on ${\mathcal X}\setminus\varSigma_{f}$.
Following \cite[Sect.\,5]{FLY08}, we recall the K\"ahler extension of $\lambda(\Omega^{p}_{{\mathcal X}/C})|_{\varDelta^{*}}$.
On ${\mathcal X}|_{{\varDelta}^{*}}$, we have the following exact sequence of holomorphic vector bundles:
$$
0
\longrightarrow
f^{*}\Omega^{1}_{C}
\longrightarrow
\Omega^{1}_{\mathcal X}
\longrightarrow
\Omega^{1}_{{\mathcal X}/C}
\longrightarrow
0.
$$
Since ${\rm rk}\,f^{*}\Omega^{1}_{C}=1$, this short exact sequence induces the following exact sequence of holomorphic vector bundles
on ${\mathcal X}|_{{\varDelta}^{*}}$:
\begin{equation}
\label{eqn:resolution:Omega:p}
0
\longrightarrow
{\mathcal E}_{{\mathcal X}/C}^{p}
\longrightarrow
\Omega^{p}_{{\mathcal X}/C}
\longrightarrow
0,
\end{equation}
where ${\mathcal E}_{{\mathcal X}/C}^{p}$ is the complex of holomorphic vector bundles over ${\mathcal X}$ given by
$$
{\mathcal E}_{{\mathcal X}/C}^{p}
\colon
(f^{*}\Omega^{1}_{C})^{\otimes p}
\to
\Omega^{1}_{\mathcal X}\otimes(f^{*}\Omega^{1}_{C})^{\otimes p-1}
\to\cdots\to
\Omega^{p-1}_{\mathcal X}\otimes f^{*}\Omega^{1}_{C}
\to
\Omega_{\mathcal X}^{p}
$$
and the map $\Omega_{\mathcal X}^{p}\to\Omega_{{\mathcal X}/C}^{p}$ is given by the canonical quotient map.
Here, if $\theta$ is a local generator of $\Omega^{1}_{C}$, then the map
$\Omega_{\mathcal X}^{i}\otimes(f^{*}\Omega^{1}_{C})^{\otimes(p-i)}\to\Omega_{\mathcal X}^{i+1}\otimes(f^{*}\Omega^{1}_{C})^{\otimes(p-i-1)}$
is given by
$\omega\otimes(f^{*}\theta)^{\otimes(p-i)}\mapsto(\omega\wedge f^{*}\theta)\otimes(f^{*}\theta)^{\otimes(p-i-1)}$
for $\omega\in\Omega_{\mathcal X}^{p-i}$.

\begin{definition}
The {\em K\"ahler extension} of $\lambda(\Omega_{{\mathcal X}/C}^{p})$ is the holomorphic line bundle over $C$ defined as
$$
\lambda({\mathcal E}_{{\mathcal X}/C}^{p})
:=
\bigotimes_{i=0}^{p}\lambda
\left(
\Omega_{\mathcal X}^{p-i}\otimes(f^{*}\Omega^{1}_{C})^{\otimes i}
\right)^{(-1)^{i}}.
$$
\end{definition}

By the exactness of \eqref{eqn:resolution:Omega:p} on ${\mathcal X}|_{{\varDelta}^{*}}$, we have the canonical isomorphism of holomorphic line bundles
over $\varDelta^{*}$:
\begin{equation}
\label{eqn:canonical:isomorphism}
\lambda(\Omega^{p}_{{\mathcal X}/C})|_{\varDelta^{*}}\cong\lambda({\mathcal E}_{{\mathcal X}/C}^{p})|_{{\varDelta}^{*}}.
\end{equation}
When $p=0,1$, the canonical isomorphism \eqref{eqn:canonical:isomorphism} extends to an isomorphism of holomorphic line bundles over $C$:
\begin{equation}
\label{eqn:canonical:isomorphism:whole:disc}
\lambda({\mathcal O}_{{\mathcal X}})=\lambda({\mathcal E}_{{\mathcal X}/C}^{0}),
\qquad
\lambda(\Omega^{1}_{{\mathcal X}/C})\cong\lambda({\mathcal E}_{{\mathcal X}/C}^{1}).
\end{equation}
\par
Via the canonical isomorphism \eqref{eqn:canonical:isomorphism}, the $L^{2}$-metric on $\lambda(\Omega^{p}_{{\mathcal X}/C})|_{\varDelta^{*}}$
induces a Hermitian metric on $\lambda({\mathcal E}_{{\mathcal X}/C}^{p})|_{{\varDelta}^{*}}$,
which is denoted by $\|\cdot\|_{\lambda({\mathcal E}_{{\mathcal X}/C}^{p}),L^{2}}$.
Notice that $\|\cdot\|_{\lambda({\mathcal E}_{{\mathcal X}/C}^{p}),L^{2}}$ does {\em not} coincide with the Hermitian metric 
on $\lambda({\mathcal E}_{{\mathcal X}/C}^{p})$ defined as the product of the $L^{2}$-metrics on 
$\lambda(\Omega_{\mathcal X}^{p-i}\otimes(f^{*}\Omega^{1}_{C})^{\otimes i})$.

\subsection
{Asymptotic behavior of the Quillen metric on $\lambda({\mathcal E}^{p}_{{\mathcal X}/C})|_{\varDelta}$}
\label{sect:3.4}
\par
Via the canonical isomorphism \eqref{eqn:canonical:isomorphism}, the line bundle $\lambda({\mathcal E}_{{\mathcal X}/C}^{p})|_{{\varDelta}^{*}}$
is endowed with the Quillen metric $\|\cdot\|_{\lambda({\mathcal E}_{{\mathcal X}/C}^{p}),Q}$.
Following \cite[Sect.\,5]{FLY08}, we recall the singularity of $\|\cdot\|_{\lambda({\mathcal E}_{{\mathcal X}/C}^{p}),Q}$ at $t=0$.
For this, we introduce some tautological vector bundles over the projective bundle ${\bf P}(T{\mathcal X})^{\lor}$.
\par
Let ${\bf P}(T{\mathcal X})^{\lor}$ be the projective bundle over ${\mathcal X}$ with projection $\varPi\colon{\bf P}(T{\mathcal X})^{\lor}\to{\mathcal X}$
such that $\varPi^{-1}(x)={\bf P}(T_{x}{\mathcal X})^{\lor}$. Here, for a complex vector space $V$, ${\bf P}(V)^{\lor}$ denotes the projective space of
hyperplanes of $V$ passing through the origin. Then we have the canonical isomorphism ${\bf P}(T{\mathcal X})^{\lor}\cong{\bf P}(\Omega_{\mathcal X}^{1})$.
We define the {\em Gauss map} $\mu\colon{\mathcal X}\setminus\varSigma_{f}\to{\bf P}(T{\mathcal X})^{\lor}$ by
$$
\mu(x):=[T_{x}X_{f(x)}],
$$
where $[T_{x}X_{f(x)}]\in{\bf P}(T_{x}{\mathcal X})^{\lor}$ is the point corresponding to the hyperplane $T_{x}X_{f(x)}\subset T_{x}{\mathcal X}$.
Then $\mu$ is a meromorphic map from $X$ to ${\bf P}(T{\mathcal X})^{\lor}$. 
Let 
$$
\sigma\colon\widetilde{\mathcal X}\to{\mathcal X}
$$ 
be a resolution of the indeterminacy of $\mu$.
Namely, there exists a birational holomorphis map
$\sigma\colon\widetilde{\mathcal X}\to{\mathcal X}$ inducing an isomorphism between $\widetilde{\mathcal X}\setminus\sigma^{-1}(\varSigma_{f})$ and
${\mathcal X}\setminus\varSigma_{f}$ such that the composite morphism
$$
\widetilde{\mu}:=\mu\circ\sigma
$$
extends to a holomorphic map from $\widetilde{\mathcal X}$ to ${\bf P}(T{\mathcal X})^{\lor}$.
We set
$$
{\rm Exc}(\sigma):=\sigma^{-1}(\varSigma_{f}).
$$
Without loss of generality, we may and will assume that ${\rm Exc}(\sigma)$ is a normal crossing divisor of $\widetilde{\mathcal X}$.

\begin{remark}
\label{remark:formula:Gauss:map}
Since $f\colon{\mathcal X}|_{\varDelta}\to\varDelta$ is a semi-stable degeneration, for any $p\in X_{0}$, there is a system of local coordinates
$(z_{0},z_{1},z_{2},z_{3})$ of ${\mathcal X}$ centered at $p$ such that
$$
f(z)=z_{0}\cdots z_{k}
\qquad
(k\leq3).
$$
Near $p$, the Gauss map $\mu$ is expressed as the following explicit meromorphic map 
$$
\mu(z)=\left(\frac{1}{z_{0}}:\cdots:\frac{1}{z_{k}}:0:\cdots:0\right).
$$
One can resolve the indeterminacy of $\mu$ in the canonical way as follows. 
\par
Let $X_{0}=E_{1}+\cdots+E_{m}$ be the irreducible decomposition. 
For $i_{1}<\cdots<i_{k}$, we set $E_{i_{1}\cdots i_{k}}:=E_{i_{1}}\cap\cdots\cap E_{i_{k}}$. 
Then $E_{i_{1}\cdots i_{k}}$ is a (possibly disconnected) submanifold of ${\mathcal X}$ and $E_{i_{1}\cdots i_{k}}=\emptyset$ for $k\geq5$.
The indeterminacy locus of $\mu$ is given by $\bigcup_{i<j}E_{ij}$. Let $\sigma^{(1)}\colon{\mathcal X}^{(1)}\to{\mathcal X}$ be the blowing-up of
$\bigcup_{i<j<k<l}E_{ijkl}$ and set $\mu^{(1)}:=\mu\circ\sigma^{(1)}$. Let $E^{(1)}_{i}$ be the proper transform of $E_{i}$ and set
$E^{(1)}_{i_{1}\cdots i_{k}}:=E^{(1)}_{i_{1}}\cap\cdots\cap E^{(1)}_{i_{k}}$ for $i_{1}<\cdots<i_{k}$.
Then $E^{(1)}_{i_{1}\cdots i_{k}}=\emptyset$ for $k\geq4$ and $E^{(1)}_{i_{1}\cdots i_{k}}$ is the proper transform of $E_{i_{1}\cdots i_{k}}$ for $k\leq3$.
The indeterminacy locus of $\mu^{(1)}$ is given by $\bigcup_{i<j}E^{(1)}_{ij}$.
Let $\sigma^{(2)}\colon{\mathcal X}^{(2)}\to{\mathcal X}^{(1)}$ be the blowing-up of $\bigcup_{i<j<k}E^{(1)}_{ijk}$ and set $\mu^{(2)}:=\mu^{(1)}\circ\sigma^{(2)}$.
Let $E^{(2)}_{i}$ be the proper transform of $E^{(1)}_{i}$ and set $E^{(2)}_{i_{1}\cdots i_{k}}:=E^{(2)}_{i_{1}}\cap\cdots\cap E^{(2)}_{i_{k}}$ for $i_{1}<\cdots<i_{k}$.
Then $E^{(2)}_{i_{1}\cdots i_{k}}=\emptyset$ for $k\geq3$ and $E^{(2)}_{i_{1}\cdots i_{k}}$ is the proper transform of $E^{(1)}_{i_{1}\cdots i_{k}}$ for $k\leq2$.
The indeterminacy locus of $\mu^{(2)}$ is given by $\bigcup_{i<j}E^{(2)}_{ij}$.
Finally, let $\sigma^{(3)}\colon{\mathcal X}^{(3)}\to{\mathcal X}^{(2)}$ be the blowing-up of $\bigcup_{i<j}E^{(2)}_{ij}$ and set $\mu^{(3)}:=\mu^{(2)}\circ\sigma^{(3)}$.
Then $\mu^{(3)}\colon{\mathcal X}^{(3)}\to{\bf P}(T{\mathcal X})^{\lor}$ is regular. Setting $\widetilde{\mathcal X}:={\mathcal X}^{(3)}$,
$\sigma:=\sigma^{(1)}\circ\sigma^{(2)}\circ\sigma^{(3)}$ and $\widetilde{\mu}:=\mu\circ\sigma$, we get a resolution of the indeterminacy of 
$\mu\colon{\mathcal X}\dashrightarrow{\bf P}(T{\mathcal X})^{\lor}$.
\end{remark}

Let $U$ be the universal hyperplane bundle over ${\bf P}(T{\mathcal X})^{\lor}$ and 
let $H$ be the universal quotient line bundle over ${\bf P}(T{\mathcal X})^{\lor}$. 
Then we have the following exact sequence of holomorphic vector bundles over ${\bf P}(T{\mathcal X})^{\lor}$:
$$
0
\longrightarrow
U
\longrightarrow
\varPi^{*}T{\mathcal X}
\longrightarrow
H
\longrightarrow
0.
$$
After \cite[Th.\,5.4]{FLY08}, we introduce the rational number $a_{p}(f,\varSigma_{f})\in{\bf Q}$ by
$$
a_{p}(f,\varSigma_{f})
:=
\sum_{j=0}^{p}(-1)^{p-j}\,
\int_{{\rm Exc}(\widetilde{\sigma})}
\widetilde{\mu}^{*}
\left\{
{\rm Td}(U)\,
\frac{{\rm Td}(c_{1}(H))-e^{-(p-j)c_{1}(H)}}{c_{1}(H)}
\right\}\,
\sigma^{*}{\rm ch}(\Omega^{j}_{\mathcal X}).
$$
Then $a_{p}(f,\varSigma_{f})$ is determined by the function germ of $f\colon{\mathcal X}|_{\varDelta}\to\varDelta$ near $\varSigma_{f}$.

\begin{theorem}
\label{thm:singularity:Quillen:metric:Kahler:extension}
Let $0\leq p\leq3$ and let $\varsigma_{p}$ be a nowhere vanishing holomorphic section of $\lambda({\mathcal E}_{{\mathcal X}/C}^{p})|_{\varDelta}$.
Then the following holds as $t\to0$:
$$
\log
\|\varsigma_{p}(t)\|^{2}_{\lambda({\mathcal E}^{p}_{{\mathcal X}/C}),Q}
=
a_{p}(f,\varSigma_{f})\,\log|t|^{2}+O(1).
$$
\end{theorem}

\begin{pf}
See \cite[Th.\,5.4]{FLY08}.
\end{pf}

\subsection
{Asymptotic behavior of the $L^{2}$-metric on $\lambda(\Omega^{p}_{{\mathcal X}/C})|_{\varDelta^{*}}$: the cases $p=2,3$}
\label{sect:3.5}
\par
Following \cite[Th.\,8.1 and Prop.\,9.6]{FLY08}, we determine the singularity of the $L^{2}$-metric on $\lambda({\mathcal E}_{{\mathcal X}/C}^{p})|_{\varDelta}$
for the remaining cases $p=2,3$.

\begin{proposition}
\label{prop:L2:metric:det:cohomology:q=2,3}
Let $\varsigma_{p}$ be a nowhere vanishing holomorphic section of $\lambda({\mathcal E}_{{\mathcal X}/C}^{p})|_{\varDelta}$.
\begin{itemize}
\item[(1)]
When $p=2$, the following holds as $t\to0$:
$$
\begin{aligned}
\log
\|\varsigma_{2}(t)\|^{2}_{\lambda({\mathcal E}^{2}_{{\mathcal X}/C}),L^{2}}
&=
\left\{
a_{2}(f,\varSigma_{f})-a_{1}(f,\varSigma_{f})+\chi(Rf_{*}{\mathcal Q}|_{\varDelta})
\right\}
\log|t|^{2}
\\
&\quad
+
O\left(\log(-\log|t|)\right).
\end{aligned}
$$
\item[(2)]
When $p=3$, the following holds as $t\to0$:
$$
\log
\|\varsigma_{3}(t)\|^{2}_{\lambda({\mathcal E}^{3}_{{\mathcal X}/C}),L^{2}}
=
\left\{
a_{3}(f,\varSigma_{f})-a_{0}(f,\varSigma_{f})
\right\}
\log|t|^{2}
+
O\left(\log(-\log|t|)\right).
$$
\end{itemize}
\end{proposition}

\begin{pf}
Let $0\leq p\leq3$. Then we have
\begin{equation}
\label{eqn:difference:L2:norm}
\begin{aligned}
\,&
\log
\left\|
\varsigma_{p}(t)
\right\|^{2}_{\lambda({\mathcal E}^{p}_{{\mathcal X}/C}),L^{2}}
-
\log
\left\|
\varsigma_{3-p}(t)
\right\|^{2}_{\lambda({\mathcal E}^{3-p}_{{\mathcal X}/C}),L^{2}}
\\
&=
\log
\left\|
\varsigma_{p}(t)\otimes\varsigma_{3-p}(t)^{\lor}
\right\|^{2}_{\lambda({\mathcal E}^{p}_{{\mathcal X}/C})\otimes\lambda({\mathcal E}^{3-p}_{{\mathcal X}/C})^{\lor},L^{2}}
\\
&=
\log
\left\|
\varsigma_{p}(t)\otimes\varsigma_{3-p}(t)^{\lor}
\right\|^{2}_{\lambda({\mathcal E}^{p}_{{\mathcal X}/C})\otimes\lambda({\mathcal E}^{3-p}_{{\mathcal X}/C})^{\lor},Q}
\\
&=
\log
\left\|
\varsigma_{p}(t)
\right\|^{2}_{\lambda({\mathcal E}^{p}_{{\mathcal X}/C}),Q}
-
\log
\left\|
\varsigma_{3-p}(t)
\right\|^{2}_{\lambda({\mathcal E}^{3-p}_{{\mathcal X}/C}),Q}
\\
&=
\left\{
a_{p}(f,\varSigma_{f})-a_{3-p}(f,\varSigma_{f})
\right\}
\log|t|^{2}
+
O(1),
\end{aligned}
\end{equation}
where the second equality follows from \cite[Eq.\,(8.4)]{FLY08} and the last equality follows from Theorem~\ref{thm:singularity:Quillen:metric:Kahler:extension}.
The result for $p=2$ follows from \eqref{eqn:difference:L2:norm} and Proposition~\ref{prop:L2:metric:det:cohomology:q=1}.
The result for $p=3$ follows from \eqref{eqn:difference:L2:norm} and Proposition~\ref{prop:L2:metric:det:cohomology:q=0}.
\end{pf}

\subsection
{Asymptotic behavior of BCOV torsion}
\label{sect:3.6}
\par
Following \cite[Th.\,8.2 and Sect.\,9.2]{FLY08}, we determine the singularity of $T_{\rm BCOV}(X_{t},g_{t})$ as $t\to0$.

\begin{theorem}
\label{thm:singularity:BCOV:torsion}
The following holds as $t\to0$:
$$
\begin{aligned}
\log T_{\rm BCOV}(X_{t},g_{t})
&=
\left\{
-3a_{0}(f,\varSigma_{f})+2a_{1}(f,\varSigma_{f})-\chi(Rf_{*}{\mathcal Q}|_{\varDelta})
\right\}
\log|t|^{2}
\\
&\qquad
+
O\left(\log(-\log|t|)\right).
\end{aligned}
$$
\end{theorem}

\begin{pf}
For simplicity, write $a_{p}$ for $a_{p}(f,\varSigma_{f})$.
Let $\varsigma_{p}$ be a nowhere vanishing holomorphic section of $\lambda({\mathcal E}_{{\mathcal X}/C}^{p})|_{\varDelta}$.
By Theorem~\ref{thm:singularity:Quillen:metric:Kahler:extension}, we get
\begin{equation}
\label{eqn:asymptotics:Quillen:BCOV:line}
\sum_{p=0}^{3}(-1)^{p}p
\log
\left\|
\varsigma_{p}(t)
\right\|^{2}_{\lambda({\mathcal E}^{p}_{{\mathcal X}/C}),Q}
=
\left(
\sum_{p=0}^{3}(-1)^{p}p\,a_{p}
\right)
\log|t|^{2}
+
O(1).
\end{equation}
By Propositions~\ref{prop:L2:metric:det:cohomology:q=0}, \ref{prop:L2:metric:det:cohomology:q=1}, \ref{prop:L2:metric:det:cohomology:q=2,3},
we get
\begin{equation}
\label{eqn:asymptotics:L2:BCOV:line}
\begin{aligned}
\sum_{p=0}^{3}(-1)^{p}p
\log
\left\|
\varsigma_{p}(t)
\right\|^{2}_{\lambda({\mathcal E}^{p}_{{\mathcal X}/C}),L^{2}}
&=
\left\{
3a_{0}-2a_{1}+2a_{2}-3a_{3}+\chi(Rf_{*}{\mathcal Q}|_{\varDelta})
\right\}
\log|t|^{2}
\\
&\qquad
+
O\left(\log(-\log|t|)\right).
\end{aligned}
\end{equation}
By the definition of Quillen metrics, we have
\begin{equation}
\label{eqn:relation:BCOV:torsion:Quillen:L2}
\log T_{\rm BCOV}(X_{t},g_{t})
=
\sum_{p=0}^{3}(-1)^{p}p
\left\{
\log
\left\|
\varsigma_{p}(t)
\right\|^{2}_{\lambda({\mathcal E}^{p}_{{\mathcal X}/C}),Q}
-
\log
\left\|
\varsigma_{p}(t)
\right\|^{2}_{\lambda({\mathcal E}^{p}_{{\mathcal X}/C}),L^{2}}
\right\}.
\end{equation}
Substituting \eqref{eqn:asymptotics:Quillen:BCOV:line}, \eqref{eqn:asymptotics:L2:BCOV:line} into \eqref{eqn:relation:BCOV:torsion:Quillen:L2}, 
we get the result.
\end{pf}

\subsection
{Asymptotic behavior of the Bott-Chern term}
\label{sect:3.7}
\par
Following \cite[Prop.\,7.9]{FLY08}, we determine the singularity of $A(X_{t},g_{t})$ as $t\to0$.
\par
Let $\varXi\in\Gamma({\mathcal X}|_{\varDelta},K_{\mathcal X})$ be a canonical form on ${\mathcal X}|_{\varDelta}$ such that (cf. \cite[Lemma 7.7]{FLY08})
\begin{equation}
\label{eqn:divisor:canonical:form}
{\rm div}(\varXi)\subset X_{0}.
\end{equation}
Let $\omega_{X_{t}}$ be the dualizing sheaf of $X_{t}$. Then $\omega_{X_{t}}\cong K_{{\mathcal X}/C}|_{X_{t}}$ for all $t\in\varDelta$.
We define $\eta_{t}\in H^{0}(X_{t},\omega_{X_{t}})$
as the canonical form on $X_{t}$ such that
$$
\varXi|_{X_{t}}=\eta_{t}\wedge df.
$$
Let $\eta_{{\mathcal X}/\varDelta}\in\Gamma(\varDelta,f_{*}K_{{\mathcal X}/C})$ be the section defined by $\eta(t)=\eta_{t}$ for all $t\in\varDelta$.
Since the family $f\colon{\mathcal X}|_{\varDelta}\to\varDelta$ is a semi-stable degeneration, $\eta_{{\mathcal X}/\varDelta}$ is regarded as a holomorphic section
of the Hodge bundle ${\mathcal F}^{3}\subset{\mathcal H}^{3}$. 
If $\varXi$ vanishes identically on $X_{0}$, then there exists $\nu\in{\bf Z}_{>0}$ such that
$t^{-\nu}\eta_{{\mathcal X}/\varDelta}$ is a nowhere vanishing holomorphic section of ${\mathcal F}^{3}$. 
Replacing $\varXi$ by $f^{*}t^{-\nu}\cdot\varXi$ in this case, we may and will assume that 
\begin{equation}
\label{eqn:non-vanishing}
\varXi|_{X_{0}}\in H^{0}(X_{0},K_{\mathcal X}|_{X_{0}})\setminus\{0\}.
\end{equation}
Namely, there is at least one irreducible component of $X_{0}$, on which $\varXi$ does not vanish.
Then $\eta_{{\mathcal X}/\varDelta}$ is a nowhere vanishing holomorphic section of ${\mathcal F}^{3}$ and hence
\begin{equation}
\label{eqn:asymptotics:relative:canonical:form}
\log\|\eta_{{\mathcal X}/\varDelta}(t)\|_{L^{2}}=O\left(\log(-\log|t|)\right)
\qquad
(t\to0)
\end{equation}
by Proposition~\ref{prop:estimate:L2:norm:determinant}.

\begin{lemma}
\label{lemma:invariance:zero:divisor}
The divisor ${\rm div}(\varXi)$ is independent of the choice of $\varXi\in\Gamma({\mathcal X}|_{\varDelta},K_{{\mathcal X}/C})$ 
satisfying \eqref{eqn:divisor:canonical:form}, \eqref{eqn:non-vanishing}.
\end{lemma}

\begin{pf}
Let $\varXi\otimes(f^{*}dt)^{-1}$ and $\varXi'\otimes(f^{*}dt)^{-1}$ be holomorphic $4$-forms on ${\mathcal X}$ satisfying 
\eqref{eqn:divisor:canonical:form}, \eqref{eqn:non-vanishing}.
Then the ratio $\varXi/\varXi'$ descends to a nowhere vanishing holomorphic function on $\varDelta^{*}$. Since both $\varXi$ and $\varXi'$ correspond
to nowhere vanishing holomorphic section of the {\em line} bundle ${\mathcal F}^{3}$, we conclude that $\varXi/\varXi'$ is a nowhere vanishing holomorphic function
on $\varDelta$. Hence ${\rm div}(\varXi)={\rm div}(\varXi')$.
\end{pf}

After Lemma~\ref{lemma:invariance:zero:divisor}, the following definition makes sense.

\begin{definition}
\label{def:normalized:canonical:divisor}
The {\em normalized canonical divisor} ${\frak K}_{({\mathcal X},X_{0})}$ of $({\mathcal X},X_{0})$ is defined as 
$$
{\frak K}_{({\mathcal X},X_{0})}:={\rm div}(\varXi),
\qquad
{\rm Supp}({\frak K}_{({\mathcal X},X_{0})})\subsetneq{\rm Supp}(X_{0}),
$$
where $\varXi$ satisfies \eqref{eqn:divisor:canonical:form}, \eqref{eqn:non-vanishing}.
\end{definition}

\par
To describe the asymptotic behavior of $A(X_{t},g_{t})$ as $t\to0$, we use the notation in Section~\ref{sect:3.4}.
Recall that $\sigma\colon\widetilde{\mathcal X}\to{\mathcal X}$ is a resolution of the indeterminacy of the Gauss map 
$\mu\colon{\mathcal X}\setminus\varSigma_{f}\to{\bf P}(T{\mathcal X})^{\lor}$ as in Remark~\ref{remark:formula:Gauss:map}, 
that $\widetilde{\mu}\colon\widetilde{\mathcal X}\to{\bf P}(T{\mathcal X})^{\lor}$ is the resolved Gauss map, 
and that $U\to{\bf P}(T{\mathcal X})^{\lor}$ is the universal hyperplane bundle. 
\par
We set 
$$
\widetilde{f}:=f\circ\sigma
$$ 
and we get a new family $\widetilde{f}\colon\widetilde{\mathcal X}\to C$, 
whose central fiber $\widetilde{X}_{0}:=\widetilde{f}^{-1}(0)$ is a possibly {\em non-reduced} normal crossing divisor. 
Hence $\widetilde{f}\colon\widetilde{\mathcal X}\to C$ is not necessarily a semi-stable degeneration.
Let $\varSigma_{\widetilde{f}}$ be the divisor of $\widetilde{\mathcal X}$ defined as the critical locus of $\widetilde{f}$:
If $\widetilde{X}_{0}=\sum_{i=1}^{k}m_{i}E_{i}$ with $E_{i}$ being an irreducible divisor of $\widetilde{\mathcal X}$ and $m_{i}\in{\bf Z}_{>0}$, then
$$
\varSigma_{\widetilde{f}}:={\rm div}(d\widetilde{f})=\sum_{i=1}^{k}(m_{i}-1)E_{i}.
$$
Since the resolution $\sigma\colon\widetilde{\mathcal X}\to{\mathcal X}$ is canonically defined, 
$\varSigma_{\widetilde{f}}$ is determined by the function germ of $f$ near $\varSigma_{f}$.

\begin{proposition}
\label{prop:asymptotics:Bott:Chern}
The following holds as $t\to0$:
$$
\log A(X_{t},g_{t})
=
-\frac{1}{12}
\left(
\int_{\sigma^{*}{\frak K}_{({\mathcal X},X_{0})}-\varSigma_{\widetilde{f}}}\widetilde{\mu}^{*}c_{3}(U)
\right)
\log|t|^{2}
+
O\left(\log(-\log|t|)\right).
$$
\end{proposition}

\begin{pf}
Let $\chi(X_{\rm gen})$ be the topological Euler number of a general fiber of $f\colon{\mathcal X}\to C$.
Let $g_{U}$ be the Hermitian metric on $U$ induced from the Hermitian metric $\varPi^{*}g^{\mathcal X}$ on $\varPi^{*}T{\mathcal X}$
via the inclusion $U\subset\varPi^{*}T{\mathcal X}$ and let $c_{3}(U)$ be the top Chern form of $(U,g_{U})$.
Define the function $A({\mathcal X}/\varDelta)$ on $\varDelta^{*}$ by $A({\mathcal X}/\varDelta)(t):=A(X_{t},g_{t})$.
By \cite[Eq.\,(7.12)]{FLY08}, we have
$$
\begin{aligned}
\log A({\mathcal X}/\varDelta)
&=
-\frac{1}{12}\widetilde{f}_{*}
\left[
\log\sigma^{*}\left(\frac{\|\varXi\|^{2}}{\|df\|^{2}}\right)\,\widetilde{\mu}^{*}c_{3}(U,g_{U})
\right]
+
\frac{\chi(X_{\rm gen})}{12}\log\|\eta_{{\mathcal X}/\varDelta}\|_{L^{2}}^{2}
\\
&=
-\frac{1}{12}\widetilde{f}_{*}
\left[
\log\left(\frac{\|\sigma^{*}\varXi\|^{2}}{\|d\widetilde{f}\|^{2}}\right)\,\widetilde{\mu}^{*}c_{3}(U,g_{U})
\right]
+
O\left(\log(-\log|t|)\right)
\\
&=
-\frac{1}{12}
\left(
\int_{\sigma^{*}{\frak K}_{({\mathcal X},X_{0})}-\varSigma_{\widetilde{f}}}\widetilde{\mu}^{*}c_{3}(U,g_{U})
\right)
\log|t|^{2}
+
O\left(\log(-\log|t|)\right),
\end{aligned}
$$
where the second equality follows from \eqref{eqn:asymptotics:relative:canonical:form}
and the third equality follows from \cite[Cor.\,4.6]{Yoshikawa07} and the equalities of divisors ${\rm div}(\sigma^{*}\varXi)=\sigma^{*}{\frak K}_{({\mathcal X},X_{0})}$,
$\varSigma_{\widetilde{f}}={\rm div}(d\widetilde{f})$ on $\widetilde{\mathcal X}$.
This completes the proof.
\end{pf}

\subsection
{Asymptotic behavior of BCOV invariants for semi-stable degenerations}
\label{sect:3.8}
\par
Define 
$$
\rho(f,\varSigma_{f})
:=
-3a_{0}(f,\varSigma_{f})+2a_{1}(f,\varSigma_{f})-\chi(Rf_{*}{\mathcal Q}|_{\varDelta})
+
\frac{1}{12}\int_{\varSigma_{\widetilde{f}}}\widetilde{\mu}^{*}c_{3}(U)
\in
{\bf Q},
$$
$$
\kappa(f,\varSigma_{f},{\frak K}_{({\mathcal X},X_{0})})
:=
\int_{\pi^{*}{\frak K}_{({\mathcal X},X_{0})}}\widetilde{\mu}^{*}c_{3}(U)
\in
{\bf Z}.
$$
Since there is a canonical way of resolving the indeterminacy of the Gauss map $\mu$ for the semi-stable degeneration
$f\colon{\mathcal X}|_{\varDelta}\to\varDelta$ as explained in Remark~\ref{remark:formula:Gauss:map},
$\rho(f,\varSigma_{f})$ is determined by the function germ of $f$ near $\varSigma_{f}$, 
whereas $\kappa(f,\varSigma_{f},{\frak K}_{({\mathcal X},X_{0})})$ is determined by the function germ of $f$ near $\varSigma_{f}$
and the normalized canonical divisor ${\frak K}_{({\mathcal X},X_{0})}$.

\begin{theorem}
\label{thm:asymptotics:BCOV:invariant:semistable:degeneration}
The following holds as $t\to0$:
$$
\log\tau_{\rm BCOV}(X_{t})
=
\left\{
\rho(f,\varSigma_{f})
-
\frac{1}{12}\kappa(f,\varSigma_{f},{\frak K}_{({\mathcal X},X_{0})})
\right\}
\log|t|^{2}
+
O\left(\log(-\log|t|)\right).
$$
\end{theorem}

\begin{pf}
Since the K\"ahler metric $g_{t}$ is induced from the K\"ahler metric $g^{\mathcal X}$ on ${\mathcal X}$, the functions on $\varDelta^{*}$
$$
t\mapsto{\rm Vol}(X_{t},g_{t}),
\qquad
t\mapsto{\rm Vol}_{L^{2}}(H^{2}(X_{t},{\bf Z}),[c_{1}({\mathcal L}_{t}])
$$
are constant by \cite[Lemma 4.12]{FLY08}. Hence there is a constant $C>0$ such that 
\begin{equation}
\label{eqn:BCOV:torsion:polarized:family}
\log\tau_{\rm BCOV}(X_{t})
=
\log T_{\rm BCOV}(X_{t},g_{t})
+
\log A(X_{t},g_{t})
+
C
\end{equation}
for all $t\in\varDelta^{*}$. 
Substituting the formulae in Theorem~\ref{thm:singularity:BCOV:torsion} and Proposition~\ref{prop:asymptotics:Bott:Chern}
into \eqref{eqn:BCOV:torsion:polarized:family}, we get the result.
\end{pf}

\subsection
{Asymptotic behavior of BCOV invariants for general degenerations}
\label{sect:3.9}
\par
By Theorem~\ref{thm:asymptotics:BCOV:invariant:semistable:degeneration}, we get the rationality of the coefficient of the logarithmic divergence
of $\log\tau_{\rm BCOV}$ for general one-parameter degenerations. In this subsection, we do not assume that $f\colon{\mathcal X}|_{\varDelta}\to\varDelta$
is a semi-stable degeneration.

\begin{theorem}
\label{thm:log:divergence:BCOV:invariant}
Let $f\colon{\mathcal X}\to C$ be a surjective morphism from an irreducible projective fourfold ${\mathcal X}$ to a compact Riemann surface $C$.
If there is a finite subset $\Delta_{f}\subset C$ such that $f|_{C\setminus\Delta_{f}}\colon{\mathcal X}|_{C\setminus\Delta_{f}}\to C\setminus\Delta_{f}$ 
is a smooth morphism and such that $X_{t}:=f^{-1}(t)$ is a Calabi-Yau threefold for all $t\in C\setminus\Delta_{f}$,
then for every $0\in\Delta_{f}$, there exists a {\em rational} number $\alpha\in{\bf Q}$ such that
$$
\log\tau_{\rm BCOV}(X_{t})=\alpha\,\log|t|^{2}+O\left(\log(-\log|t|)\right)
\qquad
(t\to0),
$$
where $t$ is a local parameter of $C$ centered at $0$. Let $g\colon({\mathcal Y},Y_{0})\to(B,0)$ be a semi-stable reduction of $f\colon({\mathcal X},X_{0})\to(C,0)$:
$$
\begin{CD}
({\mathcal Y},Y_{0})
@>\Phi>> 
({\mathcal X},X_{0})
\\
@V g VV  
@VV f V 
\\
(B,0)
@> \phi >> 
(C,0).
\end{CD}
$$ 
Then $\alpha$ is given by 
$$
\alpha
=
\frac{1}{\deg\{\phi\colon(B,0)\to(C,0)\}}
\left\{
\rho(g,\varSigma_{g})
-
\frac{1}{12}\kappa(g,\varSigma_{g},{\frak K}_{({\mathcal Y},Y_{0})})
\right\}.
$$
\end{theorem}

\begin{pf}
By the definition of semi-stable reduction \cite[Chap.\,II]{Mumford73}, ${\mathcal Y}$ is a smooth projective fourfold and $B$ is a compact Riemann surface
such that ${\mathcal Y}|_{B^{*}}\cong{\mathcal X}|_{C^{*}}\times_{C^{*}}B^{*}$ and the divisor $Y_{0}=g^{-1}(0)$ is reduced and normal crossing.
Here we set $B^{*}:=B\setminus\{0\}$ and $C^{*}:=C\setminus\{0\}$.
By choosing an appropriate local parameter $s$ of $(B,0)$, we may assume that $\phi(s)=s^{\nu}$.
Since $Y_{s}\cong X_{\phi(s)}=X_{s^{\nu}}$, the result follows from Theorem~\ref{thm:asymptotics:BCOV:invariant:semistable:degeneration} applied to 
the semi-stable degeneration of Calabi-Yau threefolds $g\colon({\mathcal Y},Y_{0})\to(B,0)$.
\end{pf}

In \cite[Th.\,9.1]{FLY08}, a weaker version of Theorem~\ref{thm:log:divergence:BCOV:invariant} was proved, where $\alpha$ was inexplicit and real.

\section
{A locality of the logarithmic singularity} 
\label{sect:4}
\par
In this section, we prove a certain locality of the cefficient $\alpha$ in Theorem~\ref{thm:log:divergence:BCOV:invariant}.
\\
\par
{\bf Set up }
Let ${\mathcal X}$ and ${\mathcal X}'$ be normal irreducible projective fourfolds. 
Let $C$ and $C'$ be compact Riemann surfaces.
Let $f\colon{\mathcal X}\to C$ and $f'\colon{\mathcal X}'\to C'$ be surjective holomorphic maps.
Let $\overline{\varSigma}_{f|_{{\mathcal X}\setminus{\rm Sing}\,{\mathcal X}}}$ 
(resp. $\overline{\varSigma}_{f'|_{{\mathcal X}'\setminus{\rm Sing}\,{\mathcal X}'}}$)
be the closure of the critical locus of $f|_{{\mathcal X}\setminus{\rm Sing}\,{\mathcal X}}$
(resp. $f'|_{{\mathcal X}'\setminus{\rm Sing}\,{\mathcal X}'}$) in ${\mathcal X}$ (resp. ${\mathcal X}'$).
Define the critical loci of $f$ and $f'$ as
$$
\varSigma_{f}:={\rm Sing}\,{\mathcal X}\cup\overline{\varSigma}_{f|_{{\mathcal X}\setminus{\rm Sing}\,{\mathcal X}}},
\qquad
\varSigma_{f'}:={\rm Sing}\,{\mathcal X}'\cup\overline{\varSigma}_{f'|_{{\mathcal X}'\setminus{\rm Sing}\,{\mathcal X}'}}
$$
and the discriminant loci of $f$ and $f'$ as
$$
\Delta_{f}:=f(\varSigma_{f}),
\qquad
\Delta_{f'}:=f'(\varSigma_{f'}).
$$
Let $0\in\Delta_{f}$ and $0'\in\Delta_{f'}$. Let $V$ (resp. $V'$) be a neighborhood of $0$ (resp. $0'$) in $C$ (resp. $C'$) such that
$V\cong\varDelta$ and $V\cap\Delta_{f}=\{0\}$ (resp. $V'\cong\varDelta$ and $V'\cap\Delta_{f'}=\{0\}$). 
In the rest of this section, we make the following:
\newline{\bf Assumption } 
\begin{itemize}
\item[(A1)]
$\Delta_{f}\not=C$, $\Delta_{f'}\not=C'$, $\dim\varSigma_{f}\leq2$, $\dim\varSigma_{f'}\leq2$, and $X_{0}$ and $X'_{0'}$ are irreducible.
\item[(A2)]
$X_{t}$ and $X'_{t'}$ are Calabi-Yau threefolds for all $t\in C\setminus\Delta_{f}$ and $t'\in C'\setminus\Delta_{f'}$.
\item[(A3)]
$f^{-1}(V)\setminus\varSigma_{f}$ carries a nowhere vanishing canonical form $\varXi$.
Similarly, $(f')^{-1}(V')\setminus\varSigma_{f'}$ carries a nowhere vanishing canonical form $\varXi'$.
\item[(A4)]
The function germ of $f$ near $\varSigma_{f}\cap f^{-1}(V)$ and the function germ of $f'$ near $\varSigma_{f'}\cap(f')^{-1}(V')$ are isomorphic.
Namely, there exist a neighborhood $O$ of $\varSigma_{f}\cap f^{-1}(V)$ in $f^{-1}(V)$,
a neighborhood $O'$ of $\varSigma_{f'}\cap(f')^{-1}(V')$ in $(f')^{-1}(V')$,
and an isomorphism $\varphi\colon O\to O'$ such that $f|_{O}=f'\circ\varphi|_{O'}$.
\end{itemize}
By (A3), (A4), the ratio $\varphi^{*}(\varXi'|_{O'})/(\varXi|_{O})$ is a nowhere vanishing holomorphic function on $O\setminus\varSigma_{f}$.
By (A1) and the normality of ${\mathcal X}$, the ratio $\varphi^{*}(\varXi'|_{O'})/(\varXi|_{O})$ extends to a nowhere vanishing holomorphic function on $O$.
Hence we have the following equality of divisors on $O$
\begin{equation}
\label{eqn:equality:normalized:canonical:divisor:1}
{\rm div}(\varXi)=\varphi^{*}{\rm div}(\varXi').
\end{equation}
\par
For $z\in C$ and $z'\in C'$, we set $X_{z}:=f^{-1}(z)$ and $X'_{z'}:=(f')^{-1}(z')$.
For $z\in C\setminus\Delta_{f}$ and $z'\in C'\setminus\Delta_{f'}$, the BCOV invariants $\tau_{\rm BCOV}(X_{z})$ and $\tau_{\rm BCOV}(X'_{z'})$ are well defined.
Let $0\in\Delta_{f}$ and $0'\in\Delta_{f'}$. A local parameter of $C$ (resp. $C'$) centered at $0$ (resp. $0'$) is denoted by $t$.
Hence $t$ is a generator of the maximal ideal of ${\mathcal O}_{C,0}$ and ${\mathcal O}_{C',0'}$.
By Theorem~\ref{thm:log:divergence:BCOV:invariant}, the functions
$t\mapsto\log\tau_{\rm BCOV}(X_{t})$ and $t\mapsto\log\tau_{\rm BCOV}(X'_{t})$ have logarithmic singularities at $0$ and $0'$, respectively.

\begin{theorem}
\label{thm:locality:BCOV:invariant}
Under {\rm (A1)--(A4)}, 
$\log\tau_{\rm BCOV}(X_{t})$ and $\log\tau_{\rm BCOV}(X'_{t})$ have the same logarithmic singularities at $t=0$:
$$
\lim_{t\to0}\frac{\log\tau_{\rm BCOV}(X_{t})}{\log|t|}
=
\lim_{t\to0}\frac{\log\tau_{\rm BCOV}(X'_{t})}{\log|t|}.
$$
In particular,
$$
\log\tau_{\rm BCOV}(X_{t})-\log\tau_{\rm BCOV}(X'_{t})=O\left(\log(-\log|t|)\right)
\qquad
(t\to0).
$$
\end{theorem}

\begin{pf}
{\em (Step 1) }
By Hironaka, there exists a succession of blowing-ups $\sigma\colon\widetilde{\mathcal X}\to{\mathcal X}$ inducing an isomorphism between
$\widetilde{\mathcal X}\setminus\sigma^{-1}(\varSigma_{f})$ and ${\mathcal X}\setminus\varSigma_{f}$ such that 
$\widetilde{X}_{0}:=(f\circ\sigma)^{-1}(0)$ is a normal crossing divisor of $\widetilde{\mathcal X}$. 
Let $\widetilde{X}_{0}=D_{0}\cup D_{1}\cup\cdots\cup D_{l}$ be the irreducible decomposition. We may and will assume that 
all $D_{\alpha}$'s are smooth hypersurfaces of $\widetilde{\mathcal X}$ and that $D_{0}$ is the proper transform of $X_{0}$.
Then $D_{1}\cup\cdots\cup D_{l}=\sigma^{-1}(\varSigma_{f})\subset\sigma^{-1}(O)$,
and $\sigma$ induces an isomorphism from $D_{0}\setminus\bigcup_{\alpha>0}D_{\alpha}$ to $X_{0}\setminus\varSigma_{f}$.
\par
Identify the pairs $(O,\varSigma_{f})$ and $(O',\varSigma_{f'})$ via $\varphi$.
We set 
$$
\widetilde{\mathcal X}':=({\mathcal X}'\setminus\varSigma_{f'})\cup\sigma^{-1}(O),
$$ 
where $\sigma^{-1}(O\setminus\varSigma_{f})$ and $O'\setminus\varSigma_{f'}$ are identified by the isomorphism $\varphi\circ\sigma$.
Then $\widetilde{\mathcal X}'$ is a smooth fourfold equipped with the projection $\sigma'\colon\widetilde{\mathcal X}'\to{\mathcal X}'$ 
defined by $\sigma':={\rm id}$ on ${\mathcal X}'\setminus\varSigma_{f'}$ and by $\varphi\circ\sigma$ on $\sigma^{-1}(O)$. 
Since $\sigma\colon\widetilde{\mathcal X}\to{\mathcal X}$ is a succession of blowing-ups, so is $\sigma'\colon\widetilde{\mathcal X}'\to{\mathcal X}'$.
We define $\widetilde{f}:=f\circ\sigma$ and $\widetilde{f}':=f'\circ\sigma'$, 
whose critical loci are denoted by $\varSigma_{\widetilde{f}}$ and $\varSigma_{\widetilde{f}'}$, respectively.
Then $\varSigma_{\widetilde{f}}\subset\sigma^{-1}(O)$ and $\varSigma_{\widetilde{f}'}\subset(\sigma')^{-1}(O')$.
\par
We set $\widetilde{\varphi}:={\rm id}_{\sigma^{-1}(O)}$. Then $\varphi$ lifts to an isomorphism $\widetilde{\varphi}\colon\sigma^{-1}(O)\cong(\sigma')^{-1}(O')$
such that $\widetilde{f}'\circ\widetilde{\varphi}=\widetilde{f}$ on $\sigma^{-1}(O)$.
Let $D'_{0}\subset\widetilde{\mathcal X}'$ be the proper transform of $X'_{0}\subset{\mathcal X}'$. 
Since $\widetilde{f}^{-1}(0)=D_{0}\cup D_{1}\cup\cdots\cup D_{l}$ and $\sigma^{-1}(\varSigma_{f})=D_{1}\cup\cdots\cup D_{l}$, 
we have the irreducible decomposition
$(\widetilde{f}')^{-1}(0')=D'_{0}\cup D'_{1}\cup\cdots\cup D'_{l}$ with $(\sigma')^{-1}(\varSigma_{f'})=D'_{1}\cup\cdots\cup D'_{l}\subset(\sigma')^{-1}(O)$, 
where we set $D'_{\alpha}:=\widetilde{\varphi}(D_{\alpha})$.
By \eqref{eqn:equality:normalized:canonical:divisor:1} and the equality $\sigma'\circ\widetilde{\varphi}=\varphi\circ\sigma$, we get
\begin{equation}
\label{eqn:equality:normalized:canonical:divisor:2}
{\rm div}(\sigma^{*}\varXi)=\widetilde{\varphi}^{*}{\rm div}\left((\sigma')^{*}\varXi'\right)\subset\sigma^{-1}(O).
\end{equation}
\par{\em (Step 2) }
Let $d\in{\bf Z}_{>0}$. Let $\pi_{d}\colon(C_{d},0_{d})\to(C,0)$ (resp. $\pi'_{d}\colon(C'_{d},0'_{d})\to(C',0')$) be a ramified covering with ramification index $d$
at $0_{d}\in C_{d}$ (resp. $0'_{d}\in C'_{d}$).
Let $\widetilde{\mathcal X}_{d}$ (resp. $\widetilde{\mathcal X}'_{d}$) be the normalization of the fibered product
$\widetilde{\mathcal X}\times_{C}C_{d}$ (resp. $\widetilde{\mathcal X}'\times_{C'}C'_{d}$)
and set $\widetilde{f}_{d}:={\rm pr}_{2}\colon\widetilde{\mathcal X}_{d}\to C_{d}$ 
(resp. $\widetilde{f}'_{d}:={\rm pr}_{2}\colon\widetilde{\mathcal X}'_{d}\to C'_{d}$).
Let $\widetilde{O}_{d}$ (resp. $\widetilde{O}'_{d}$) be the open subset of $\widetilde{\mathcal X}_{d}$ (resp. $\widetilde{\mathcal X}'_{d}$) defined as
${\rm pr}_{1}^{-1}(\sigma^{-1}(O))$ (resp. $({\rm pr}_{1}')^{-1}((\sigma')^{-1}(O'))$).
Then $\widetilde{\varphi}\colon\sigma^{-1}(O)\cong(\sigma')^{-1}(O')$ lifts to an isomorphism $\widetilde{\varphi}_{d}\colon\widetilde{O}_{d}\cong\widetilde{O}'_{d}$ 
such that $\widetilde{f}_{d}=\widetilde{f}'_{d}\circ\widetilde{\varphi}_{d}$.
\par
Define $U_{d}:=\widetilde{\mathcal X}_{d}\setminus\widetilde{f}_{d}^{-1}(0_{d})$ and 
$U'_{d}:=\widetilde{\mathcal X}'_{d}\setminus(\widetilde{f}')_{d}^{-1}(0'_{d})$.
By \cite[Chap.\,II, \S3]{Mumford73}, the pairs $(\widetilde{\mathcal X}_{d},U_{d})$ and $(\widetilde{\mathcal X}'_{d},U'_{d})$ are toroidal embeddings.
(See \cite[Chap.\,II \S1]{Mumford73} for the notion of toroidal embeddings.)
Let $\widetilde{f}_{d}^{-1}(0_{d})=E_{0}\cup E_{1}\cup\cdots\cup E_{m}$ 
(resp. $(\widetilde{f}'_{d})^{-1}(0'_{d})=E'_{0}\cup E'_{1}\cup\cdots\cup E'_{m}$) be the irreducible decomposition,
where $E_{0}$ (resp. $E'_{0}$) is the component corresponding to the proper transforms of $X_{0}$ (resp. $X'_{0}$) in $\widetilde{\mathcal X}$ 
(resp. $\widetilde{\mathcal X}'$). Then $m\geq n$.
Since $X_{0}\setminus\varSigma_{f}\cong D_{0}\setminus\bigcup_{\alpha>0}D_{\alpha}$ is reduced and smooth,
we have $E_{0}\setminus\bigcup_{\beta>0}E_{\beta}\cong D_{0}\setminus\bigcup_{\alpha>0}D_{\alpha}\cong X_{0}\setminus\varSigma_{f}$. 
When $\beta>0$, we deduce from \cite[Chap.\,II \S3]{Mumford73} that ${\rm pr}_{1}(E_{\beta})=D_{\alpha(\beta)}$ for some $\alpha(\beta)>0$.
Similarly, $E'_{0}\setminus\bigcup_{\beta>0}E'_{\beta}\cong X'_{0}\setminus\varSigma_{f}$ 
and ${\rm pr}_{1}(E'_{\beta})=D'_{\alpha(\beta)}$ $(\alpha(\beta)>0)$ for $\beta>0$.
\par{\em (Step 3) }
By an appropriate choice of $d\in{\bf Z}_{>0}$, there exists a sheaf of ideals ${\mathcal I}_{d}\subset{\mathcal O}_{\widetilde{\mathcal X}_{d}}$
with ${\mathcal I}_{d}|_{U_{d}}={\mathcal O}_{U_{d}}$, 
whose blowing-up $\varpi\colon{\mathcal Y}_{d}\to\widetilde{\mathcal X}_{d}$ provides a semi-stable reduction of $f\colon({\mathcal X},X_{0})\to(C,0)$, 
i.e., the following commutative diagram (cf. \cite[Chap.\,II]{Mumford73})
$$
\begin{CD}
{\mathcal Y}_{d}
@>\varpi>> 
\widetilde{\mathcal X}_{d}
@>{\rm pr}_{1}>> 
\widetilde{\mathcal X}
@>\sigma>>
{\mathcal X} 
\\
@V g_{d} VV  
@V \widetilde{f}_{d} VV 
@V \widetilde{f} VV
@V f VV
\\
C_{d} 
@> {\rm id} >> 
C_{d} 
@> \pi_{d} >> 
C
@> {\rm id} >> 
C,
\end{CD}
$$
where ${\mathcal Y}_{d}$ is smooth and $g_{d}^{-1}(0_{d})$ is a {\em reduced}, normal crossing divisor of ${\mathcal Y}_{d}$. 
We define $E_{\beta}^{o}:=E_{\beta}\setminus\bigcup_{\beta'\not=\beta}E_{\beta'}$. 
Since the ideal sheaf ${\mathcal I}_{d}$ is of the form as in \cite[p.\,91 last line]{Mumford73},
we deduce from \cite[Th.\,9*]{Mumford73} that there exists $\nu_{0}\in{\bf Z}$ with
\begin{equation}
\label{eqn:ideal:sheaf:blowup}
{\mathcal I}_{d}|_{U_{d}\cup E_{0}^{o}}\cong{\mathcal O}_{U_{d}\cup E_{0}^{o}}(\nu_{0}E_{0}^{o}).
\end{equation}
Since $E_{0}^{o}$ is a smooth divisor of $\widetilde{\mathcal X}_{d}\setminus\bigcup_{\beta>0}E_{\beta}=U_{d}\cup E_{0}^{o}$, 
we deduce from \eqref{eqn:ideal:sheaf:blowup} and the definition of blowing-up of sheaf of ideals that the maps
$\varpi\colon\varpi^{-1}(E_{0}^{o})\to E_{0}^{o}$ and
$\varpi\colon{\mathcal Y}_{d}\setminus\varpi^{-1}(E_{1}\cup\cdots\cup E_{m})\to\widetilde{\mathcal X}_{d}\setminus(E_{1}\cup\cdots\cup E_{m})$ 
are isomorphisms. 
Write $g_{d}^{-1}(0_{d})=F_{0}+\cdots+F_{n}$, where every $F_{i}$ is irreducible and $F_{0}$ is the proper transform of $E_{0}$. 
Then $n\geq m$ and $F_{1}\cup\cdots\cup F_{n}=\varpi^{-1}(E_{1}\cup\cdots\cup E_{m})\subset\varpi^{-1}(\widetilde{O}_{d})$.
\par{\em (Step 4) }
We define the sheaf of ideals ${\mathcal I}'_{d}\subset{\mathcal O}_{\widetilde{\mathcal X}'_{d}}$ by
$$
{\mathcal I}'_{d}|_{\widetilde{\mathcal X}'_{d}\setminus\bigcup_{\beta>0}E'_{\beta}}
=
{\mathcal O}_{\widetilde{\mathcal X}'_{d}\setminus\bigcup_{\beta>0}E'_{\beta}}(\nu_{0}E'_{0}),
\qquad
{\mathcal I}'_{d}|_{O'_{d}}=(\varphi_{d})_{*}{\mathcal I}_{d}.
$$
Then ${\mathcal I}'_{d}|_{U'_{d}}={\mathcal O}_{U'_{d}}$.
Let $\varpi'\colon{\mathcal Y}'\to\widetilde{\mathcal X}'_{d}$ be the blowing-up of ${\mathcal I}'_{d}$ and set $g'_{d}:=\widetilde{f}'_{d}\circ\varpi'$.
Since the map $\varpi'\colon(\varpi')^{-1}(\widetilde{O}'_{d})\to\widetilde{O}'_{d}$ is identified with 
$\varpi\colon\varpi^{-1}(\widetilde{O}_{d})\to\widetilde{O}_{d}$ via the identification $\widetilde{\varphi}_{d}\colon \widetilde{O}_{d}\cong\widetilde{O}'_{d}$, 
we get the isomorphism of divisors
$$
(g'_{d})^{-1}(0'_{d})\cap(\varpi')^{-1}(\widetilde{O}'_{d})\cong(F_{0}\cap\widetilde{O}_{d})+F_{1}+\cdots+F_{n},
$$
which is a reduced normal crossing divisor of $\varpi^{-1}(\widetilde{O}_{d})$.
Let $F'_{0}\subset(g'_{d})^{-1}(0_{d})$ be the proper transform of $E'_{0}\subset\widetilde{\mathcal X}'_{d}$ and 
let $F'_{\gamma}\subset(g'_{d})^{-1}(0_{d})$ be the irreducible component corresponding to $F_{\gamma}$ for $\gamma>0$.
Then $F'_{0}\cap(\varpi')^{-1}(\widetilde{O}'_{d})\cong F_{0}\cap\widetilde{O}_{d}$.
Since the map $\varpi'\colon{\mathcal Y}'\setminus(\varpi')^{-1}(\widetilde{O}'_{d})\to\widetilde{\mathcal X}'_{d}\setminus\widetilde{O}'_{d}$ is an isomorphism,
$F'_{0}\setminus(\varpi')^{-1}(\widetilde{O}'_{d})\cong X'_{0}\setminus\widetilde{O}'$ is a smooth divisor of ${\mathcal Y}'\setminus(\varpi')^{-1}(\widetilde{O}'_{d})$.
Hence $(g'_{d})^{-1}(0_{d})=F'_{0}+\cdots+F'_{n}$ is a reduced normal crossing divisor of ${\mathcal Y}'$.
Thus $g'_{d}\colon{\mathcal Y}'_{d}\to C_{d}$ is a semi-stable reduction of $f'\colon{\mathcal X}'\to C'$:
We have the following commutative diagram
$$
\begin{CD}
{\mathcal Y}'_{d}
@>\varpi'>> 
\widetilde{\mathcal X}'_{d}
@>{\rm pr}_{1}>> 
\widetilde{\mathcal X}'
@>\sigma'>>
{\mathcal X} '
\\
@V g'_{d} VV  
@V \widetilde{f}'_{d} VV 
@V \widetilde{f}' VV
@V f' VV
\\
C'_{d} 
@> {\rm id} >> 
C'_{d} 
@> \pi_{d} >> 
C'
@> {\rm id} >> 
C'.
\end{CD}
$$
\par{\em (Step 5) }
Set $O_{d}:=\varpi^{-1}(\widetilde{O}_{d})$ and $O'_{d}:=(\varpi')^{-1}(\widetilde{O}'_{d})$. 
Let $\varphi_{d}\colon O_{d}\cong O'_{d}$ be the isomorphism induced by $\widetilde{\varphi}_{d}\colon \widetilde{O}_{d}\cong\widetilde{O}'_{d}$.
Since
$$
\varSigma_{g_{d}}\subset O_{d},
\qquad
\varSigma_{g'_{d}}=\varphi_{d}(\varSigma_{g_{d}})\subset O'_{d},
\qquad
(g_{d},O_{d})=(g'_{d}\circ\varphi_{d},O_{d})
$$
by construction, we get the following equality by the definition in Section~\ref{sect:3.8}
\begin{equation}
\label{eqn:equality:lambda}
\rho(g_{d},\varSigma_{g_{d}})=\rho(g'_{d},\varSigma_{g'_{d}}).
\end{equation}
\par{\em (Step 6) }
Set $\psi:=\sigma\circ{\rm pr}_{1}\circ\varpi\colon{\mathcal Y}_{d}\to{\mathcal X}$ and 
$\psi':=\sigma'\circ{\rm pr}_{1}\circ\varpi'\colon{\mathcal Y}'_{d}\to{\mathcal X}'$.
Let $\Upsilon$ (resp $\Upsilon'$) be a canonical form defined near $g_{d}^{-1}(0_{d})$ (resp. $(g'_{d})^{-1}(0'_{d})$)
and satisfying \eqref{eqn:divisor:canonical:form}, \eqref{eqn:non-vanishing}. 
Then there exist $a_{\gamma},a'_{\gamma}\in{\bf Z}_{\geq0}$ for $0\leq\gamma\leq n$ such that
\begin{equation}
\label{eqn:zero:divisor:Upsilon}
{\rm div}(\Upsilon)=\sum_{\gamma=0}^{n}a_{\gamma}F_{\gamma},
\qquad
{\rm div}(\Upsilon')=\sum_{\gamma=0}^{n}a'_{\gamma}F'_{\gamma}.
\end{equation}
\par
Since $\psi^{*}\varXi$ (resp. $(\psi')^{*}\varXi'$) is a (possibly meromorphic) $4$-form defined
on a neighborhood of $g_{d}^{-1}(0_{d})$ (resp. $(g'_{d})^{-1}(0'_{d})$), whose possible zeros and poles are supported on $g_{d}^{-1}(0_{d})$ 
(resp. $(g'_{d})^{-1}(0'_{d})$), we can express
\begin{equation}
\label{eqn:zero:divisor:pullback:Xi}
{\rm div}\left(\psi^{*}\varXi\right)=\sum_{\gamma=0}^{n}b_{\gamma}F_{\gamma},
\qquad
{\rm div}\left((\psi')^{*}\varXi'\right)=\sum_{\gamma=0}^{n}b'_{\gamma}F'_{\gamma},
\end{equation}
where $b_{\gamma},b'_{\gamma}\in{\bf Z}$ for $0\leq\gamma\leq n$. 
Since $\varXi$ (resp. $\varXi'$) is nowhere vanishing on $f^{-1}(V)\setminus O$ (resp. $(f')^{-1}(V')\setminus O'$) by assumption and 
since $\psi$ (resp. $\psi'$) has ramification index $d$ along $F_{0}\setminus\bigcup_{\gamma>0}F_{\gamma}$ 
(resp. $F'_{0}\setminus\bigcup_{\gamma>0}F'_{\gamma}$), 
$\psi^{*}\varXi$ (resp. $(\psi')^{*}\varXi'$) has zeros of order $d-1$ on the proper transform of $E_{0}$ (resp. $E'_{0}$). 
Hence
\begin{equation}
\label{eqn:order:F0:F'0}
b_{0}=b'_{0}=d-1.
\end{equation}
Since the map $\varpi'\colon O'_{d}\to\widetilde{O}'_{d}$ is identified with $\varpi\colon O_{d}\to\widetilde{O}_{d}$ 
via the identification $\widetilde{\varphi}_{d}\colon \widetilde{O}_{d}\cong\widetilde{O}'_{d}$, 
we get by \eqref{eqn:equality:normalized:canonical:divisor:2}
$$
{\rm div}(\psi^{*}\varXi)\cap O_{d}=\varphi_{d}^{*}({\rm div}((\psi')^{*}\varXi'))\cap O_{d}.
$$
Hence
\begin{equation}
\label{eqn:order:Fi:F'i}
b_{\gamma}=b'_{\gamma}
\qquad
(\gamma>0).
\end{equation}
\par
Since $\varphi^{*}(\varXi'\wedge\overline{\varXi}')/(\varXi\wedge\overline{\varXi})$ is a nowhere vanishing positive function on $O$
by \eqref{eqn:equality:normalized:canonical:divisor:1},
there exist constants $C_{1},C_{2}>0$ such that for all $t$ with $0<|t|\ll1$,
\begin{equation}
\label{eqn:estimate:integral:relative:canonical:form:1}
C_{1}\int_{X_{t}\cap O}\sqrt{-1}\frac{\varXi\wedge\overline{\varXi}}{df\wedge d\overline{f}}
\leq
\int_{X'_{t}\cap O'}\sqrt{-1}\frac{\varXi'\wedge\overline{\varXi}'}{df'\wedge d\overline{f}'}
\leq
C_{2}\int_{X_{t}\cap O}\sqrt{-1}\frac{\varXi\wedge\overline{\varXi}}{df\wedge d\overline{f}}.
\end{equation}
Since $\sqrt{-1}\varXi\wedge\overline{\varXi}$ (resp. $\sqrt{-1}\varXi'\wedge\overline{\varXi}'$) is a volume form on $f^{-1}(V)\setminus O$
(resp. $(f')^{-1}(V')\setminus O'$) and since $f$ (resp. $f'$) has no critical points on $f^{-1}(V)\setminus O$ (resp. $(f')^{-1}(V')\setminus O'$),
there exist constants $C_{3},C_{4}>0$ such that for all $t$ with $0<|t|\ll1$,
\begin{equation}
\label{eqn:estimate:integral:relative:canonical:form:2}
C_{3}\int_{X_{t}\setminus O}\sqrt{-1}\frac{\varXi\wedge\overline{\varXi}}{df\wedge d\overline{f}}
\leq
\int_{X'_{t}\setminus O'}\sqrt{-1}\frac{\varXi'\wedge\overline{\varXi}'}{df'\wedge d\overline{f}'}
\leq
C_{4}\int_{X_{t}\setminus O}\sqrt{-1}\frac{\varXi\wedge\overline{\varXi}}{df\wedge d\overline{f}}.
\end{equation}
By \eqref{eqn:estimate:integral:relative:canonical:form:1}, \eqref{eqn:estimate:integral:relative:canonical:form:2}, 
there exist constants $C_{5},C_{6}>0$ such that for all $t$ with $0<|t|\ll1$,
\begin{equation}
\label{eqn:estimate:integral:relative:canonical:form:3}
C_{5}
\left\|
\left.\frac{\varXi}{df}\right|_{X_{t}}
\right\|_{L^{2}}^{2}
\leq
\left\|
\left.\frac{\varXi'}{df'}\right|_{X'_{t}}
\right\|_{L^{2}}^{2}
\leq
C_{6}
\left\|
\left.\frac{\varXi}{df}\right|_{X_{t}}
\right\|_{L^{2}}^{2}.
\end{equation}
\par{\em (Step 7) }
Let $s$ be a local parameter of $C_{d}$ and $C'_{d}$ centered at $0_{d}$ and $0'_{d}$.
Since $\pi_{d}\colon C_{d}\to C$ and $\pi'_{d}\colon C'_{d}\to C'$ has ramification index $d$ at $s=0$, we may assume 
\begin{equation}
\label{eqn:relation:local:parameters}
\pi_{d}(s)=s^{d}.
\end{equation}
By the definition of $\Upsilon$ (resp. $\Upsilon'$), the map $s\mapsto(\Upsilon/dg_{d})|_{Y_{s}}\in H^{0}(Y_{s},K_{Y_{s}})$
(resp. $s\mapsto(\Upsilon'/dg'_{d})|_{Y'_{s}}\in H^{0}(Y'_{s},K_{Y'_{s}})$)
is a nowhere vanishing holomorphic section of $(g_{d})_{*}K_{{\mathcal Y}_{d}/C_{d}}$ (resp. $(g'_{d})_{*}K_{{\mathcal Y}'_{d}/C'_{d}}$) near $0_{d}$ (resp. $0'_{d}$).
By \eqref{eqn:asymptotics:relative:canonical:form}, we get 
\begin{equation}
\label{eqn:estimate:integral:relative:canonical:form:5}
\log\left\|
\left.\frac{\Upsilon}{dg_{d}}\right|_{Y_{s}}
\right\|_{L^{2}}^{2}
=
O\left(\log(-\log|s|)\right),
\qquad
\log\left\|
\left.\frac{\Upsilon'}{dg'_{d}}\right|_{Y'_{s}}
\right\|_{L^{2}}^{2}
=
O\left(\log(-\log|s|)\right)
\end{equation}
as $s\to 0$, where we set $Y_{s}:=g_{d}^{-1}(s)$ and $Y'_{s}:=(g'_{d})^{-1}(s)$.
\par
Since the fibers $Y_{s}$ and $Y'_{s}$ are Calabi-Yau threefolds for $s\not=0$,
the map $s\mapsto\psi^{*}\varXi/dg_{d}|_{Y_{s}}$ (resp. $s\mapsto(\psi')^{*}\varXi'/dg'_{d}|_{Y'_{s}}$) is a holomorphic section of
$(g_{d})_{*}K_{{\mathcal Y}_{d}/C_{d}}$ (resp. $(g'_{d})_{*}K_{{\mathcal Y}'_{d}/C'_{d}}$) near $0_{d}$ (resp. $0'_{d}$). 
Hence there exist $c,c'\in{\bf Z}$ and $\epsilon(s),\epsilon'(s)\in{\mathcal O}(\varDelta)$ such that
\begin{equation}
\label{eqn:ratio:canonical:forms}
\left.\frac{\psi^{*}\varXi}{\Upsilon}\right|_{Y_{s}}=s^{c}\epsilon(s),
\qquad
\left.\frac{(\psi')^{*}\varXi'}{\Upsilon'}\right|_{Y'_{s}}=s^{c'}\epsilon'(s),
\qquad
\epsilon(0)\not=0,
\quad
\epsilon'(0)\not=0.
\end{equation}
By \eqref{eqn:zero:divisor:Upsilon}, \eqref{eqn:zero:divisor:pullback:Xi}, \eqref{eqn:order:F0:F'0}, \eqref{eqn:order:Fi:F'i}, \eqref{eqn:ratio:canonical:forms}, we get
\begin{equation}
\label{eqn:estimate:integral:relative:canonical:form:4}
a_{\gamma}=b_{\gamma}+c,
\qquad
a'_{\gamma}=b'_{\gamma}+c'=b_{\gamma}+c'.
\end{equation}
Since $Y_{s}=X_{\pi_{d}(s)}=X_{s^{d}}$ and $Y'_{s}=X'_{s^{d}}$ for $s\not=0$ by \eqref{eqn:relation:local:parameters}, we get
by \eqref{eqn:estimate:integral:relative:canonical:form:5}, \eqref{eqn:ratio:canonical:forms},
\begin{equation}
\label{eqn:estimate:integral:relative:canonical:form:6}
\begin{aligned}
\log\left\|
\left.\frac{\varXi}{df}\right|_{X_{s^{d}}}
\right\|_{L^{2}}^{2}
&=
\log\left\|
\psi^{*}\left\{\left.\frac{\varXi}{df}\right|_{X_{s^{d}}}\right\}
\right\|_{L^{2}}^{2}
=
\log\left\|
\left.\frac{\varXi}{d(g_{d}^{d})}\right|_{Y_{s}}
\right\|_{L^{2}}^{2}
\\
&=
\log
\left(
|s|^{-2(d-1)}
\left\|
\left.\frac{\varXi}{dg_{d}}\right|_{Y_{s}}
\right\|_{L^{2}}^{2}
\right)
+
O(1)
\\
&=
-(d-1)\,\log|s|^{2}
+
\log
\left|
\left.\frac{\psi^{*}\varXi}{\Upsilon}\right|_{Y_{s}}
\right|^{2}
+
\log\left\|
\left.\frac{\Upsilon}{dg_{d}}\right|_{Y_{s}}
\right\|_{L^{2}}^{2}
+
O(1)
\\
&=
(c-d+1)\log|s|^{2}+O\left(\log(-\log|s|)\right)
\qquad
(s\to0).
\end{aligned}
\end{equation}
Similarly, we get
\begin{equation}
\label{eqn:estimate:integral:relative:canonical:form:7}
\log\left\|
\left.\frac{\varXi'}{df'}\right|_{X'_{s^{d}}}
\right\|_{L^{2}}^{2}
=
(c'-d+1)\log|s|^{2}+O\left(\log(-\log|s|)\right)
\qquad
(s\to0).
\end{equation}
Comparing \eqref{eqn:estimate:integral:relative:canonical:form:3} and
\eqref{eqn:estimate:integral:relative:canonical:form:6}, \eqref{eqn:estimate:integral:relative:canonical:form:7}, we get
\begin{equation}
\label{eqn:vanishing:order:relative:canonical:form}
c=c'.
\end{equation}
By \eqref{eqn:estimate:integral:relative:canonical:form:4}, \eqref{eqn:vanishing:order:relative:canonical:form}, 
we get $a_{\gamma}=a'_{\gamma}$ for all $0\leq\gamma\leq n$. 
Hence
\begin{equation}
\label{eqn:comparison:normalized:canonical:divisor:1}
{\rm div}(\Upsilon)=\sum_{\gamma=0}^{n}a_{\gamma}F_{\gamma},
\qquad
{\rm div}(\Upsilon')=\sum_{\gamma=0}^{n}a_{\gamma}F'_{\gamma},
\qquad
a_{\gamma}\in{\bf Z}_{\geq0},
\quad
(\forall\,\gamma\geq0).
\end{equation}
Since $F_{0}\cap O_{d}$ (resp. $F_{\gamma}$ ($\gamma>0$)) is identified with $F'_{0}\cap O'_{d}$ (resp. $F'_{\gamma}$ ($\gamma>0$)) via $\varphi_{d}$,
we get by \eqref{eqn:comparison:normalized:canonical:divisor:1} and the definition of normalized canonical divisor in Section~\ref{sect:3.8}
the following equality of divisors via the identification $\varphi_{d}\colon O_{d}\cong O'_{d}$:
\begin{equation}
\label{eqn:comparison:normalized:canonical:divisor:2}
{\frak K}_{({\mathcal Y}_{d},Y_{0})}\cap O_{d}=\varphi_{d}^{*}({\frak K}_{({\mathcal Y}'_{d},Y'_{0})}\cap O'_{d}).
\end{equation}
Since $(g_{d},\varSigma_{g_{d}})=(g'_{d}\circ\varphi_{d},\varSigma_{g_{d}})$ and $\varSigma_{d}=\varphi_{d}(\varSigma_{g'_{d}})$ 
via $\varphi_{d}\colon O_{d}\cong O'_{d}$, 
we deduce from \eqref{eqn:comparison:normalized:canonical:divisor:2} and the definition in Section~\ref{sect:3.8} the equality
\begin{equation}
\label{eqn:equality:kappa}
\kappa(g_{d},\varSigma_{g_{d}},{\frak K}_{({\mathcal Y}_{d},Y_{0})})
=
\kappa(g'_{d},\varSigma_{g'_{d}},{\frak K}_{({\mathcal Y}'_{d},Y'_{0})}).
\end{equation}
By Theorem~\ref{thm:asymptotics:BCOV:invariant:semistable:degeneration} and \eqref{eqn:equality:lambda}, \eqref{eqn:equality:kappa},
we get
\begin{equation}
\label{eqn:log:divergence:semistable:reduction}
\lim_{s\to 0}\frac{\log\tau_{\rm BCOV}(Y_{s})}{\log|s|^{2}}
=
\lim_{s\to 0}\frac{\log\tau_{\rm BCOV}(Y'_{s})}{\log|s|^{2}}.
\end{equation}
Since $Y_{s}=X_{s^{d}}$ and $Y'_{s}=X'_{s^{d}}$ for $s\not=0$, the result follows from \eqref{eqn:log:divergence:semistable:reduction}.
\end{pf}

\section
{Degenerations to Calabi-Yau varieties with ordinary double points} 
\label{sect:5}
\par
In this section, we determine the asymptotic behavior of BCOV invariants for the simplest degenerations of Calabi-Yau threefolds, i.e.,
degenerations to Calabi-Yau varieties with at most ordinary double points (cf. \cite[\S2]{FLY08}).
Recall that an $n$-dimensional singularity is an {\it ordinary double point} if it is isomorphic to the hypersurface singularity 
at $0\in{\bf C}^{n}$ defined by the equation $z_{0}^{2}+\cdots+z_{n}^{2}=0$.

\begin{definition}
A complex projective variety $X$ of dimension $3$ is a {\em Calabi-Yau variety with at most ordinary double points} if the following are satisfied: 
\begin{itemize}
\item[(1)]
There exists a nowhere vanishing canonical form on $X\setminus{\rm Sing}(X)$.
\item[(2)]
$X$ is connected and $H^{q}(X,{\mathcal O}_{X})=0$ for $0<q<3$.
\item[(3)]
${\rm Sing}(X)$ consists of at most ordinary double points.
\end{itemize}
\end{definition}

\begin{theorem}
\label{thm:singularity:BCOV:invariant:ODP}
Let $f\colon{\mathcal X}\to C$ be a surjective morphism from a smooth projective fourfold ${\mathcal X}$ to a compact Riemann surface $C$.
Let $\Delta_{f}\subset C$ be the discriminant locus of $f\colon{\mathcal X}\to C$ and assume that $X_{t}:=f^{-1}(t)$ is a Calabi-Yau threefold
for all $t\in C\setminus\Delta_{f}$.
Let $0\in\Delta_{f}$ and let $t$ be a local parameter of $C$ centered at $0$.
If $X_{0}$ is a Calabi-Yau variety with at most ordinary double points, then
$$
\log\tau_{\rm BCOV}(X_{t})=\frac{\#{\rm Sing}\,X_{0}}{6}\,\log|t|^{2}+O\left(\log(-\log|t|)\right)
\qquad
(t\to0).
$$
\end{theorem}

\begin{pf}
{\em (Step 1) }
Since the deformation germ $f\colon({\mathcal X},X_{0})\to(C,0)$ is a smoothing of $X_{0}$, we have $h^{1,2}(X_{0})=h^{2,1}(X_{0})\geq 1$. 
Since $t\mapsto h^{2}(X_{t},\Omega^{1}_{X_{t}})$ is a locally constant function on $C$ by \cite[Th.\,2.11]{FLY08}, we get
$h^{1,2}(X_{t})=h^{2,1}(X_{t})\geq 1$.
When $h^{1,2}(X_{t})=h^{2,1}(X_{t})=1$ and $\#{\rm Sing}\,X_{0}=1$, the result was proved in \cite[Th.\,8.2]{FLY08}.
Since there does exit a family of Calabi-Yau threefolds $f'\colon{\mathcal X}'\to C'$ with $h^{1,2}(X'_{t})=1$ such that ${\rm Sing}\,X'_{0}$ consists 
of a unique ordinary double point, e.g. the family of quintic mirror threefolds (cf. \cite[\S12]{FLY08}),
we get the result by Theorem~\ref{thm:locality:BCOV:invariant} and \cite[Th.\,8.2]{FLY08} when $\#{\rm Sing}\,X_{0}=1$.
\par{\em (Step 2) }
Fix a family of Calabi-Yau threefolds over a compact Riemann surface ${\frak f}\colon{\frak X}\to{\frak C}$ 
such that ${\frak X}_{0}$, $0\in{\frak C}$, is a Calabi-Yau variety with a unique ordinary double point as its singular set. 
Fix its semi-stable reduction ${\frak g}\colon({\frak Y},{\frak Y}_{0})\to({\frak B},0)$: We have a commutative diagram:
$$
\begin{CD}
({\frak Y},{\frak Y}_{0})
@>\Psi>> 
({\frak X},{\frak X}_{0})
\\
@V {\frak g} VV  
@VV {\frak f} V 
\\
({\frak B},0)
@> \psi >> 
({\frak C},0).
\end{CD}
$$
Let $\nu$ be the ramification index of $\psi\colon({\frak B},0)\to({\frak C},0)$.
Let ${\frak Y}_{0}={\frak E}_{0}+\cdots+{\frak E}_{n}$ be the irreducible decomposition
such that ${\frak E}_{0}\setminus\bigcup_{\alpha>0}{\frak E}_{\alpha}\cong{\frak X}_{0}\setminus{\rm Sing}\,{\frak X}_{0}$. 
Then $\Psi$ ramifies along ${\frak E}_{0}$ with ramification index $\nu$.
There exist a neighborhood ${\frak U}$ of ${\rm Sing}\,{\frak X}_{0}$ in ${\frak X}$ and an open subset ${\frak V}$ of ${\frak Y}$ such that
${\frak V}=\Psi^{-1}({\frak U})$ and ${\frak E}_{1}\cup\cdots\cup{\frak E}_{n}\subset{\frak V}$.
Then $\Psi$ induces an isomorphism between ${\frak E}_{0}\setminus{\frak V}$ and ${\frak X}_{0}\setminus{\frak U}$
and the map $\Psi$ has ramification index $\nu$ on ${\frak E}_{0}\setminus\bigcup_{\alpha>0}{\frak E}_{\alpha}$.
\par{\em (Step 3) }
By an appropriate choices of local parameters $t$ of $({\frak C},0)$ and $s$ of $({\frak B},0)$, we may assume $t=s^{\nu}$. 
By identifying ${\frak f}$ with $t\circ{\frak f}$ and ${\frak g}$ with $s\circ{\frak g}$, we have the following equality of functions defined near ${\frak Y}_{0}$:
$$
\Psi^{*}{\frak f}={\frak g}^{\nu}.
$$
\par
Since ${\frak X}_{0}$ is a Calabi-Yau variety with a unique ordinary double point as its singular set and since ${\frak X}$ is smooth,
there exists a nowhere vanishing holomorphic $4$-form $\varXi$ defined on a neighborhood of ${\frak X}_{0}$ in ${\frak X}$.
Since ${\frak X}_{0}$ has only canonical singularities, the function 
$$
t\mapsto
\left\|(\varXi/d{\frak f})|_{{\frak X}_{t}}\right\|_{L^{2}}^{2}
$$
is continuous around $0\in{\frak C}$ and $\left\|(\varXi/d{\frak f})|_{{\frak X}_{0}}\right\|_{L^{2}}^{2}\not=0$ by e.g. \cite[Th.\,7.2]{Yoshikawa10}.
Since
$$
\nu\,\Psi^{*}(\varXi/d{\frak f})=\Psi^{*}\varXi/({\frak g}^{\nu-1}d{\frak g}),
$$
the fact that the function $t\mapsto\left\|(\varXi/d{\frak f})|_{{\frak X}_{t}}\right\|_{L^{2}}^{2}$ is $C^{0}$ and does not vanish at $t=0$ implies that the map
$s\mapsto(\Psi^{*}\varXi/{\frak g}^{\nu-1}d{\frak g})|_{{\frak Y}_{s}}\in H^{0}({\frak Y}_{s},K_{\frak Y}|_{{\frak Y}_{s}})$ is a nowhere vanishing holomorphic section of
${\frak g}_{*}K_{\frak Y}$ defined near $s=0$. Hence $(\Psi^{*}\varXi)/{\frak g}^{\nu-1}$ is a nowhere vanishing holomorphic $4$-form defined 
near ${\frak Y}_{0}$ and 
$$
0
\leq
{\frak K}_{({\frak Y},{\frak Y}_{0})}
=
{\rm div}(\Psi^{*}\varXi/{\frak g}^{\nu-1})={\rm div}(\Psi^{*}\varXi)-(\nu-1){\frak Y}_{0}.
$$
Write ${\rm div}(\Psi^{*}\varXi)=\sum_{\alpha=0}^{n}a_{\alpha}\,{\frak E}_{\alpha}$, $a_{\alpha}\in{\bf Z}_{\geq0}$.
Since ${\frak Y}_{0}={\rm div}({\frak g})=\sum_{\alpha=0}^{n}{\frak E}_{\alpha}$ by the reducedness of ${\frak Y}_{0}$, we get by the effectivity of
${\rm div}(\Psi^{*}\varXi)-(\nu-1){\frak Y}_{0}$
$$
a_{\alpha}\geq\nu-1
\qquad
(\forall\,\alpha\in\{0,1,\ldots,n\}).
$$
On the other hand, since $\Psi$ has ramification index $\nu$ on ${\frak E}_{0}\setminus\bigcup_{\alpha>0}{\frak E}_{\alpha}$, we have $\alpha_{0}=\nu-1$.
Hence we can express
\begin{equation}
\label{eqn:normalized:canonical:divisor:canonical:singularity:1}
{\frak K}_{({\frak Y},{\frak Y}_{0})}
=
{\rm div}(\Psi^{*}\varXi/{\frak g}^{\nu-1}|_{\frak V})
=
\sum_{\alpha=1}^{n}b_{\alpha}{\frak E}_{\alpha}
\qquad
(b_{\alpha}\in{\bf Z}_{\geq0}).
\end{equation}
\par{\em (Step 4) }
Consider the general case. Set $m:=\#{\rm Sing}\,X_{0}\geq1$. Then ${\rm Sing}\,X_{0}=\{p_{1},\ldots,p_{m}\}$. 
Since every germ $f\in{\mathcal O}_{X_{0},p_{i}}$ is isomorphic to the germ $(z_{0})^{2}+(z_{1})^{2}+(z_{2})^{2}+(z_{3})^{2}$ at $0\in{\bf C}^{4}$,
it follows from the construction of semi-stable reduction \cite[Chap.\,II \S3]{Mumford73} that
there exists a semi-stable reduction 
$$
\begin{CD}
({\mathcal Y},Y_{0})
@>\Phi>> 
({\mathcal X},X_{0})
\\
@V g VV  
@VV f V 
\\
(B,0)
@> \phi >> 
(C,0)
\end{CD}
$$
with the following properties:
\begin{itemize}
\item[(i)]
The ramification index of $\phi\colon(B,0)\to(C,0)$ is given by $\nu$.
\item[(ii)]
For every $p_{i}\in{\rm Sing}\,X_{0}$, there is a neighborhood $U_{i}$ of $p_{i}$ in ${\mathcal X}$ such that
$(f,U_{i})\cong({\frak f},{\frak U})$ and $(g,\Phi^{-1}(U_{i}))\cong({\frak g},{\frak V})$. 
\item[(iii)]
Set $V_{i}:=\Phi^{-1}(U_{i})$. Then $Y_{0}\setminus\bigcup_{i=1}^{m}V_{i}\cong X_{0}\setminus\bigcup_{i=1}^{m}U_{i}$.
\end{itemize}
By (ii), the irreducible component of $Y_{0}=g^{-1}(0)$ contained in $V_{i}$ can be expressed as $E_{1}^{(i)}+\cdots+E_{n}^{(i)}$
and satisfy $E_{0}\cap V_{i}+E_{1}^{(i)}+\cdots+E_{n}^{(i)}\cong{\frak E}_{0}\cap{\frak V}+{\frak E}_{1}+\cdots+{\frak E}_{n}$,
where $E_{\alpha}^{(i)}\cong{\frak E}_{\alpha}$ for all $i=1,\ldots,m$ and $1\leq\alpha\leq n$.
This implies that
\begin{equation}
\label{eqn:multiplicativity:lambda}
\rho(g,\varSigma_{g})
=
\sum_{i=1}^{m}\rho(g|_{V_{i}},\varSigma_{g|_{V_{i}}})
=
m\,\rho({\frak g},\varSigma_{\frak g}).
\end{equation}
\par{\em (Step 5) }
Let $t$ be a local parameter of $(C,0)$ and let $s$ be a local parameter of $(B,0)$. As in Step 3, we may assume $\phi^{*}t=s^{\nu}$ and hence
$$
\Phi^{*}f=g^{\nu}.
$$
Let $\omega$ be a nowhere vanishing holomorphic $4$-form defined near $X_{0}$. Since $X_{0}$ has only canonical singularities,
the section of the Hodge bundle $t\mapsto\omega/df|_{X_{t}}\in H^{0}(X_{t},K_{X_{t}})$ is holomorphic and nowhere vanishing.
By the same reason as in Step 3, the section of Hodge bundle $s\mapsto\Phi^{*}\omega/(g^{\nu-1}dg)|_{Y_{s}}\in H^{0}(Y_{s},K_{Y_{s}})$
is holomorphic and nowhere vanishing around $s=0$ and we get
\begin{equation}
\label{eqn:normalized:canonical:divisor:canonical:singularity:2}
{\frak K}_{({\mathcal Y},Y_{0})}={\rm div}\left(\Phi^{*}\omega/g^{\nu-1}\right),
\qquad
{\rm Supp}\,{\frak K}_{({\mathcal Y},Y_{0})}
\subset
\bigcup_{i=1}^{m}\bigcup_{\alpha=1}^{n}E_{\alpha}^{(i)}.
\end{equation}
We may assume that, under the identification $(f,U_{i})\cong({\frak f},{\frak U})$ in (ii), 
the ratio $\varXi/\omega$ is a nowhere vanishing holomorphic function on ${\frak U}$.
Then we get by \eqref{eqn:normalized:canonical:divisor:canonical:singularity:1}, \eqref{eqn:normalized:canonical:divisor:canonical:singularity:2}
\begin{equation}
\label{eqn:normalized:canonical:divisor:canonical:singularity:3}
{\frak K}_{({\mathcal Y},Y_{0})}=\sum_{i=1}^{m}\sum_{\alpha=1}^{n}b_{\alpha}\,E_{\alpha}^{(i)}.
\end{equation}
By (ii), \eqref{eqn:normalized:canonical:divisor:canonical:singularity:3} and the definition of $\kappa(g,\varSigma_{g},{\frak K}_{({\mathcal Y},Y_{0})})$,
we get
\begin{equation}
\label{eqn:multiplicativity:kappa}
\kappa(g,\varSigma_{g},{\frak K}_{({\mathcal Y},Y_{0})})
=
m\,\kappa({\frak g},\varSigma_{\frak g},{\frak K}_{({\frak Y},{\frak Y}_{0})}).
\end{equation}
Since ${\frak X}_{s^{\nu}}={\frak Y}_{s}$, we get by Theorem~\ref{thm:asymptotics:BCOV:invariant:semistable:degeneration} and Step 1
$$
\rho({\frak g},\varSigma_{\frak g})-\frac{1}{12}\kappa({\frak g},\varSigma_{\frak g},{\frak K}_{({\frak Y},{\frak Y}_{0})})
=
\lim_{s\to0}\frac{\log\tau_{\rm BCOV}({\frak Y}_{s})}{\log|s|^{2}}
=
\lim_{s\to0}\frac{\log\tau_{\rm BCOV}({\frak X}_{s^{\nu}})}{\log|s|^{2}}
=
\frac{\nu}{12}.
$$
Thus we get by \eqref{eqn:multiplicativity:lambda}, \eqref{eqn:multiplicativity:kappa}
$$
\rho(g,\varSigma_{g})-\frac{1}{12}\kappa(g,\varSigma_{g},{\frak K}_{({\mathcal Y},Y_{0})})
=
m\,\{\rho({\frak g},\varSigma_{\frak g})-\frac{1}{12}\kappa({\frak g},\varSigma_{\frak g},{\frak K}_{({\frak Y},{\frak Y}_{0})})\}
=
\frac{m\nu}{12}.
$$
By Theorem~\ref{thm:asymptotics:BCOV:invariant:semistable:degeneration} again and the relation $Y_{s}=X_{s^{\nu}}$, we get
$$
\begin{aligned}
\log\tau_{\rm BCOV}(X_{s^{d}})
&=
\log\tau_{\rm BCOV}(Y_{s})
\\
&=
\{\rho(g,\varSigma_{g})-\frac{1}{12}\kappa(g,\varSigma_{g},{\frak K}_{({\mathcal Y},Y_{0})})\}\log|s|^{2}
+
O\left(\log(-\log|s|)\right)
\\
&=
\frac{m\nu}{12}\log|s|^{2}+O\left(\log(-\log|s|)\right)
\\
&=
\frac{m}{12}\log|s^{\nu}|^{2}+O\left(\log(-\log|s|)\right)
\qquad
(s\to0).
\end{aligned}
$$
This completes the proof.
\end{pf}


\end{document}